\documentclass[12pt,a4paper]{article}
\usepackage[latin1]{inputenc}
\usepackage[T1]{fontenc}
\usepackage{amsmath}
\usepackage{amssymb}
\usepackage{textcomp}
\usepackage[english]{babel}
\usepackage[all]{xy}

\newtheorem{soulem}{Sublemma}
\newtheorem{thm}{Theorem}
\newtheorem{lem}{Lemma}
\newtheorem{prop}{Proposition}
\newtheorem{cor}{Corollary}

\newcommand\dem{\textbf{Proof: }}
\newcommand\Def{\textbf{Definition: }}
\newcommand\Defs{\textbf{Definitions: }}
\newcommand\rem{\textbf{Remark: }}

\title{Reducibility of quasiperiodic cocycles in linear Lie groups}
\author{C.Chavaudret}
\date{}

\begin{document}

\maketitle

Abstract: Let $G$ be a linear Lie group. We define the $G$-reducibility of a continuous or discrete cocycle modulo 
$N$. We show that a $G$-valued continuous or discrete cocycle which is $GL(n,\mathbb{C})$-reducible is in fact 
$G$-reducible modulo 2 if $G=GL(n,\mathbb{R}),SL(n,\mathbb{R}), Sp(n,\mathbb{R})$ or $O(n)$ 
and modulo 1 if $G=U(n)$.


\section*{Introduction}

\noindent Let $G$ be a Lie subgroup of $GL(n,\mathbb{C})$ and $\mathcal{G}$ its Lie algebra.
Let $\mathbb{T}^d=\mathbb{R}^d/\mathbb{Z}^d$ and $N\mathbb{T}^d=\mathbb{R}^d/(N\mathbb{Z})^d$ for 
$N\in\mathbb{N}\setminus\{0\}$.
Let us consider the equation

\begin{equation}\label{gen}\forall t\in\mathbb{R},\forall \theta\in\mathbb{T}^d, \frac{d}{dt}X(t,\theta)=
A(\theta+t\omega)X(t,\theta)\end{equation}

\noindent where $A:\mathbb{T}^d\rightarrow \mathcal{G}$ is continuous and $\omega\in\mathbb{R}^d$ is rationally independant. 
Let ${X:(t,\theta)\mapsto X^t(\theta)}$ be the associated continuous cocycle, 
i.e the map from $\mathbb{R}\times \mathbb{T}^d$ to $G$ satisfying (\ref{gen}) such that 
for all $\theta\in\mathbb{T}^d$, $X^0(\theta)=Id$. 
The terminology comes from the fact that $X$ satisfies the cocycle relation 

\begin{equation}\forall t,s\in\mathbb{R},\forall \theta\in\mathbb{T}^d, \ 
X^{t+s}(\theta)=X^t(\theta+s\omega)X^s(\theta)\end{equation}

\noindent As $X$ is continuous in the variable $t$, $X^t(\theta)$ remains in the connected component of the identity for all  
$t,\theta$, so we can suppose $G$ is connected.

\bigskip
\Def 
Let $X$ be a $G$-valued continuous cocycle. 
We say that $X$ is $G$-reducible modulo 
$N\in\mathbb{N}\setminus\{0\}$ 
if there exists $Z:N\mathbb{T}^d\rightarrow G$ continuous and $B\in \mathcal{G}$ such that for all 
$t\in\mathbb{R}$ and 
$\theta\in\mathbb{T}^d$,

\begin{equation}\label{ladef}X^t(\theta)=Z(\theta+t\omega)^{-1}e^{tB}Z(\theta)\end{equation}

\noindent We say $X$ is reducible if it is reducible modulo 1.

\bigskip
\rem For a continuous cocycle, reducibility implies that for all $\theta$,

\begin{equation}\partial_\omega Z(\theta)=BZ(\theta)-Z(\theta)A(\theta)\end{equation}

\noindent where $\partial_\omega Z(\theta):=\frac{d}{dt}Z(\theta+t\omega)_{|_{t=0}}$.

\bigskip
\noindent We shall prove the following theorems for continuous cocycles before adapting them to discrete cocycles:

\begin{thm}\label{CR'} Let $X$ be a continuous cocycle with values in $GL(n,\mathbb{R})$; 
if it is $GL(n,\mathbb{C})$-reducible, then it is $GL(n,\mathbb{R})$-reducible
 modulo $2$.\end{thm}

\begin{thm}\label{Sp'etc} Let $X$ be a $G$-valued continuous cocycle, where 
$G$ is either the symplectic group $Sp(n,\mathbb{R})$\footnote{In this case, $n$ is even}, the group $SL(n,\mathbb{R})$ of matrices with determinant 1, 
the orthogonal group $O(n)$, or the unitary group $U(n)$.
Suppose $X$ is
 $GL(n,\mathbb{C})$-reducible. 
Then it is $G$-reducible modulo 2 if $G=Sp(n,\mathbb{R}),
SL(n,\mathbb{R})$ or $O(n)$ and modulo 1 if $G=U(n)$. 

 \end{thm}

\bigskip
\Def Assume $(\omega,1)$ is rationally independant. A discrete $G$-valued cocycle is a map 
$X:\mathbb{Z}\times \mathbb{T}^d\rightarrow G$ 
such that for all $n,m\in\mathbb{Z}$ and all $\theta\in\mathbb{T}^d$,

\begin{equation}X^{n+m}(\theta)=X^n(\theta+m\omega)X^m(\theta)\end{equation}

\bigskip
\Def A discrete cocycle $X$ is $G$-reducible modulo $N$ if there exists a continuous 
$Z:N\mathbb{T}^d\rightarrow G$ and $A\in {G}$ 
such that 

$$\forall n\in\mathbb{Z},\forall \theta\in\mathbb{T}^d,\ X^n(\theta)=Z(\theta+n\omega)^{-1}A^nZ(\theta)$$ 

\noindent This is equivalent to the fact that $X^1(\theta)=Z(\theta+\omega)^{-1}AZ(\theta)$ for all $\theta$. 
A discrete cocycle is reducible if it is reducible modulo 1.

\bigskip

\noindent Theorems \ref{CR'} and \ref{Sp'etc} also hold for discrete cocycles. Adapting their proofs to the discrete case, 
one gets:

\begin{thm}\label{reddisc}Let $X$  
be a $G$-valued discrete cocycle with $G$ in $GL(n,\mathbb{R})$,$SL(n,\mathbb{R})$, $Sp(n,\mathbb{R}), O(n)$ or 
$U(n)$, 
and assume it is $GL(n,\mathbb{C})$-reducible. Then 
$X$ is $G$-reducible modulo $\chi_G$ with

$$\chi_G=\left\{\begin{array}{cc}
2 & \mathrm{if}\ G=GL(n,\mathbb{R}),SL(n,\mathbb{R}), G=Sp(n,\mathbb{R})\ 
\mathrm{or}\ G=O(n)\\
1 & \mathrm{if}\ G=U(n)\\
\end{array}\right.$$

\end{thm}

\bigskip
In \cite{HY}, H.He and J.You have solved a conjecture from \cite{E3} showing that if $\omega$ is diophantine, 
if $X_\lambda$ is the cocycle which is solution of 

\begin{equation}\label{perturbation}
\frac{d}{dt}X^t(\theta,\lambda)=(A(\lambda)+F_\epsilon(\theta,\lambda))X^t(\theta,\lambda)\end{equation}

\noindent where $F_\epsilon$ is sufficiently small and $A(\lambda)$ satisfies non-degeneracy conditions on an interval 
$\Lambda\subset \mathbb{R}$, then $X_\lambda$ is $GL(n,\mathbb{C})$-reducible 
for almost all $\lambda\in\Lambda$. 

\noindent Applying theorems \ref{CR'} and \ref{Sp'etc} to this result, we get that if $X^t(\theta,\lambda)$ is $G$-valued, 
with $G$ in $GL(n,\mathbb{R})$,$SL(n,\mathbb{R})$, $Sp(n,\mathbb{R}), O(n),U(n)$, 
then for almost all $\lambda\in\Lambda$,
$X^t(\theta,\lambda)$ is $G$-reducible modulo 2 if $G=GL(n,\mathbb{R}), Sp(n,\mathbb{R}), SL(n,\mathbb{R}), 
O(n)$ and modulo 1 if $G=U(n)$.

\noindent This completes R.Krikorian's result (see \cite{Kr2}): let $A(\lambda)$ be a generic one-parameter family 
taking its values in the Lie algebra of a compact semi-simple group $G$; then the system (\ref{perturbation}) 
is $G$-reducible for almost every $\lambda$ modulo some integer $\chi_G$ depending only on $G$, and $\chi_G=1$ if 
$G=U(n)$. Now we know that $\chi_G=2$ if $G=O(n)$.

\bigskip
So, when $G$ is real, there is a loss of periodicity. In the periodic case ($d=1$), this is a well-known phenomenon. 
However, it seems that there exists a large class of real cocycles 
that are reducible in a subgroup of 
$GL(n,\mathbb{R})$ without loss of periodicity. 
For instance, in \cite{Kr} (proposition 2.2.4), R. Krikorian has showed when $G$ is a compact semi-simple group
 that if a discrete cocycle $X$ is $G$-reducible modulo $m$ to a constant cocycle $n\mapsto e^{nB}$, 
then there exists a subset $S\subset G$ of Haar measure 1 such that if $e^B\in S$, then $X$ is reducible 
modulo 1. This tells us that loss of periodicity is quite rare, at least in the compact case. 

\noindent We shall prove the following: 

\begin{prop}\label{NR} If a continuous $G$-valued cocycle $X$, with $G=GL(n,\mathbb{R}),SL(n,\mathbb{R}), 
Sp(n,\mathbb{R})$ or $O(n)$, is $GL(n,\mathbb{C})$-reducible to a cocycle $t\mapsto e^{tB}$ such that no eigenvalue of 
$B$ is in $\mathbb{R}+i\pi\langle\mathbb{Z}^d,\omega\rangle\setminus\{0\}$, then $X$ is $G$-reducible.\end{prop}

\noindent There is a natural question: considering a generic one-parameter family of cocycles 
which are solution of (\ref{perturbation}), where $A(\lambda)$ satisfies non-degeneracy conditions, 
is it true that for almost all $\lambda$ 
such that the cocycle $X_\lambda$ is reducible to a constant cocycle $t\mapsto e^{tB_\lambda}$, 
no eigenvalue of $B_\lambda$ is in 
$\mathbb{R}+i\pi\langle\mathbb{Z}^d,\omega\rangle\setminus\{0\}$?

\noindent If it were true, the already mentioned result of \cite{HY} and proposition \ref{NR} would 
imply $G$-reducibility almost everywhere, without loss of periodicity, for generic one-parameter families of cocycles 
of type (\ref{perturbation}). 

\bigskip
\rem All the results which we shall prove also hold in higher regularity classes: defining "$C^r$-reducibility" in the same way as reducibility, 
but with $Z$ in $C^r$ and not only continuous, it is easy to check that 
we get theorems \ref{CR'}, \ref{Sp'etc}, \ref{reddisc} and 
proposition \ref{NR} with "$C^r$-reducibility" instead of "reducibility".

\bigskip
\textbf{Sketch of the proof}

\bigskip
We shall define notions of invariant subbundle and of Jordan subbundle as families parametrized by $\mathbb{T}^d$ 
and with values in the subspaces of $\mathbb{C}^n$, satisfying a continuity condition and some invariance properties. 
In order to prove theorem \ref{CR'}, we shall first study the properties of the decomposition of $\mathbb{C}^n$ 
into Jordan subbundles given by the $GL(n,\mathbb{C})$-reducibility of $X$ to a cocycle $t\mapsto e^{tB}$; 
we shall decompose 
$\mathbb{R}^n$ into two reducible invariant subbundles, one of them, say $W$, modulo 2 and having a basis with real exponents, 
the other, say $W'$, modulo 1 and having a basis  with no exponent in $\mathbb{R}+i\pi\langle\mathbb{Z}^d,\omega\rangle$, 
and such that the gap between the imaginary parts of two exponents cannot be in $2\pi\langle\mathbb{Z}^d,\omega\rangle$ (this is 
called a non-resonance condition). This gives theorem \ref{CR'} as a corollary, but it also facilitates the proof 
of theorem \ref{Sp'etc} for the orthogonal and the symplectic group, since it is then easy to construct 
real global bases for the cocycle's invariant subbundles, which are respectively orthonormal and 
symplectic.

\bigskip
If in equation (\ref{ladef}), no eigenvalue of $B$ has its imaginary part in $\pi
\langle\mathbb{Z}^d,{\omega}\rangle
\setminus \{0\}$, then the first of these two subbundles, $W$, is trivial, so we can have real reducibility 
without doubling the period, and consequently, 
we can also have $G$-reducibility without loss of periodicity, if $G=SL(n,\mathbb{R}),Sp(n,\mathbb{R})$ or $O(n)$, 
whence proposition \ref{NR}.

\bigskip
In order to get theorem \ref{Sp'etc} for $SL(n,\mathbb{R})$, we will just apply theorem 
\ref{CR'}, then show that the determinant of 
$Z$ is constant, so we can assume it equal to 1. Notice that no condition on the exponents of the subbundles 
is used.

\bigskip
In the case where $G=U(n)$, we shall start from the decomposition of $\mathbb{C}^n$ into complex Jordan subbundles 
with non-resonant exponents, and construct a global complex orthonormal basis. 
As $U(n)$ is not a real Lie group, we do not need to double the period.

\bigskip
To prove theorem \ref{reddisc}, we can make exactly the same proof as for theorems \ref{CR'} and 
\ref{Sp'etc}, simply adapting the first lemma to the discrete case, i.e considering integer translations 
instead of continuous translations in the direction of $\omega$. The dynamics are not modified by the fact that 
the time is discrete. 

\bigskip
For a particular class of discrete cocycles, there is another way of proving $G$-reducibility:

\bigskip
\Def A discrete cocycle $X$ is called $G$-exponential if there exists a continuous $A:\mathbb{T}^d\rightarrow \mathcal{G}$ 
such that for all $\theta\in\mathbb{T}^d$, $X^1(\theta)=e^{A(\theta)}$.

\bigskip
To prove theorem \ref{reddisc} for $G$-exponential cocycles, we can also construct a suspension of $X$ on a torus 
of greater dimension, 
taking its values in $G$, using the function $A$ from the definition of a 
$G$-exponential cocycle. We will obtain a continuous cocycle over $(\omega,1)$, which is possible since $(\omega,1)$ 
is assumed to be rationally independant. We then show that if $X$ is $GL(n,\mathbb{R})$-reducible, 
then so is its suspension. 
Using theorems \ref{CR'} and \ref{Sp'etc}, we obtain $G$-reducibility for the suspension 
modulo 1 or 2. Restricting to integer time and to a subtorus, we finally obtain $G$-reducibility for $X$ 
modulo 1 or 2.

\bigskip
\noindent
\textbf{Notations}

\bigskip
For a vector $v\in \mathbb{R}^n$, denote by 
$Re\ v$ and $Im\ v$ its real and imaginary parts. The euclidean scalar product is denoted by 
$\langle,\rangle$ for a real vector space, and $\langle,\rangle_\mathbb{C}$ for a complex vector space (we shall take it semilinear in the 
second variable); euclidean distance is denoted 
by $d(,)$. Also, we shall write $M^*$ for the adjoint of a matrix $M$; 
if $M$ is real, $M^*$ is simply the transpose of $M$. Matrix $\left(\begin{array}{cc}
0 & -I_n\\
I_n & 0\\
\end{array}\right)$ is denoted by $J$. Finally, $\mathbb{N}^*=\mathbb{N}\setminus\{0\}$.

\section{$GL(n,\mathbb{R})$-reducibility}\label{redCR}

\noindent In this section, we shall assume that $X$ is a real cocycle. 

\subsection{Preliminary lemmas}

\begin{lem}\label{dep}\begin{enumerate}
\item Let $\omega\in\mathbb{R}^d$ rationally independant, $\beta\in\mathbb{R}$ and $N\in\mathbb{N}^*$. 
Suppose that for any real sequence $t_j\rightarrow\infty$,

\begin{equation}t_j\omega\rightarrow 0\in N\mathbb{T}^d\Rightarrow t_j\beta\rightarrow 0\in 2\pi\mathbb{T}\end{equation}

\noindent Then there exists $k\in\mathbb{Z}^d$ such that $\beta=2\pi\langle k,\frac{\omega}{N}\rangle$.

\item Let $\omega\in\mathbb{R}^d$ such that $(\omega,1)$ is rationally independant and 
$N\in\mathbb{N}^*$. 
Suppose that for any integer sequence $t_j\rightarrow\infty$,

\begin{equation}t_j\omega\rightarrow 0\in N\mathbb{T}^d\Rightarrow t_j\beta\rightarrow 0\in 2\pi\mathbb{T}\end{equation}

\noindent Then there exists $k\in\mathbb{Z}^d$ such that $\beta=2\pi\langle k,\frac{(\omega,1)}{N}\rangle$.

\end{enumerate}
\end{lem}

\dem 1. It is enough to prove the assertion for $N=1$: if it is true for $N=1$ and that for every sequence 
$t_j$, $t_j\omega\rightarrow 0\in N\mathbb{T}^d\Rightarrow t_j\beta\rightarrow 0\in 2\pi\mathbb{T}$, 
then $t_j\frac{\omega}{N}\rightarrow 0\in \mathbb{T}^d\Rightarrow t_j\beta \rightarrow 0\in 2\pi\mathbb{T}$, 
and applying the case $N=1$ with $\frac{\omega}{N}$ instead of $\omega$, we get $\beta=2\pi\langle k,\frac{\omega}
{N}\rangle$ 
for some $k\in\mathbb{Z}^d$. 

\bigskip
First notice that $(\frac{\beta}{2\pi},{\omega})$ is rationally dependant, since the orbit of the translation 
in the direction $(\frac{\beta}{2\pi},{\omega})$ 
is not dense in $\mathbb{T}^{d+1}$: otherwise, there would exist a sequence $(t_j)$ satisfying 
$t_j(\frac{\beta}{2\pi},{\omega})
\rightarrow (\frac{1}{2},0)\in \mathbb{T}^{d+1}$, which would contradict the assumption. 
So there exists $k=(k_1,\dots k_d)\in\mathbb{Z}^d, p\in\mathbb{Z}$ such that $(p,k)$ is primitive (i.e the greatest common 
divisor of $k_i$ and $p$ is 1) and

\begin{equation}\label{ratdep}\langle k,{\omega}\rangle+p\frac{\beta}{2\pi}=0 \end{equation}

\noindent Notice that this is the only possible resonance (i.e $(p,k)$ is unique up to a scalar). 
For if there existed a
$(p',k')$ independant from $(p,k)$ and such that $\langle k',{\omega}\rangle+p'\frac{\beta}{2\pi}=0$, then 

\begin{equation}\left\{\begin{array}{c}
p\beta+2\pi\langle k,{\omega}\rangle=0\\
p'\beta+2\pi\langle k',{\omega}\rangle=0\\
\end{array}\right.\end{equation}

\noindent would hold, so $pk'-p'k=0\in\mathbb{R}^d$ since $\omega$ is rationally independant, 
which would contradict the assumption that $(p,k)$ and $(p',k')$ are independant.

\bigskip
Now let us show that $p=\pm 1$. By contraposition, suppose $|p|\geq 2$.

\noindent Let $V$ the subspace of $\mathbb{R}^{d+1}$ generated by $(p,k)$ .
Let $m_1,\dots m_d\in \mathbb{Z}^{d+1},m_i=(m_{i,1},\dots m_{i,d+1})$ such that the $(d+1)\times (d+1)$-matrix

\begin{equation}C:=\left(\begin{array}{c}
(p,k)\\
m_1\\
\vdots\\
m_d\\
\end{array}\right)\end{equation}

\noindent has determinant 1. Such a matrix exists, according to \cite{Ca}, corollary 3 p.14. 
Form the following commuting diagram:

\begin{equation}
\xymatrix{\relax
\mathbb{R}^{d+1}  \ar[d]_{\Pi} \ar[r]^C & \mathbb{R}^{d+1} \ar[d]_{\Pi}\\
\mathbb{T}^{d+1} \ar[r]^{\bar{C}} & \mathbb{T}^{d+1}
}\end{equation}

\noindent where $\Pi$ is the canonical projection from $\mathbb{R}^{d+1}$ onto $\mathbb{T}^{d+1}$. 
As $C$ has determinant 1, $\bar{C}$ is a homeomorphism. 
So the orbit of $\Pi\left(\begin{array}{c}
\frac{\beta}{2\pi}\\
^t{\omega}
\end{array}\right)$ is dense is $\Pi(V^\perp)$ if the orbit of 
$\Pi(C\left(\begin{array}{c}
\frac{\beta}{2\pi}\\
^t{\omega}\\
\end{array}\right))$ is dense in $\Pi(C(V^\perp))$. Now $\Pi(C(V^\perp))\subset\{0\}\times \mathbb{T}^d$ and 

$$\Pi(C\left(\begin{array}{c}
\frac{\beta}{2\pi}\\
^t{\omega}\\
\end{array}\right))= \left(\begin{array}{c}
0\\
\langle m_1,(\frac{\beta}{2\pi},{\omega})\rangle\\
\vdots \\
\langle m_d,(\frac{\beta}{2\pi},{\omega})\rangle
\end{array}\right)$$

\noindent Moreover, assume that

\begin{equation}\sum_ia_i\langle m_i,(\frac{\beta}{2\pi},{\omega})\rangle=0\Leftrightarrow
\sum_i\sum_{j\neq 1}a_im_{i,j}{\omega_{j-1}}+\sum_ia_im_{i,1}\frac{\beta}{2\pi}=0\end{equation}

\noindent then, as the resonance is unique, $(\sum_ia_im_{i,1},\dots \sum_ia_im_{i,d+1})=\sum_ia_im_i$ 
is a multiple of
 $(p,k)$, which is impossible since $m_i$ are independant from $(p,k)$ by definition. So 
 $(\langle m_{1},(\frac{\beta}{2\pi},{\omega})\rangle,\dots \langle m_{d},(\frac{\beta}{2\pi},{\omega})\rangle)$ 
 is rationally independant 
and its orbit is dense in $\mathbb{T}^{d}$. Therefore, the orbit of $\Pi\left(\begin{array}{c}
\frac{\beta}{2\pi}\\
^t{\omega}\end{array}\right)$ is 
dense in 
$\Pi(V^\perp)$. 

\bigskip
Let $m\in\mathbb{Z}^d$ such that $\frac{\langle k,m\rangle}{p}$ is not an integer (it exists, since $(p,k)$ is primitive and 
$|p|\geq 2$). 
Then $^t(\frac{\langle k,m\rangle}{p},-m)\in V^\perp$. 
As $\Pi\left(\begin{array}{c}
\frac{\beta}{2\pi}\\
^t{\omega}\end{array}\right)$ 
dans $\Pi(V^\perp)$ has a dense orbit, 
there exists an unbounded sequence $t_j$ such that $\Pi\left(\begin{array}{c}
t_j\frac{\beta}{2\pi}\\
^t{t_j\omega}\end{array}\right)\rightarrow 
\Pi\left(\begin{array}{c}
\frac{\langle k,m\rangle}{p}\\
-^tm\end{array}\right)=:\left(\begin{array}{c}
\alpha\\
0\end{array}\right)\neq \left(\begin{array}{c}
0\\
0\end{array}\right)$, which contradicts our assumption. 

\bigskip
2. Again, we can assume that $N=1$. Let $t_j$ be a real unbounded sequence such that  
$t_j(\omega,1)\rightarrow 0\in\mathbb{T}^{d+1}$. For all $j$, let $n_j\in\mathbb{Z}$ and $r_j\in [0,1[$ such that 
$t_j=n_j+r_j$. In particular, $t_j\rightarrow 0\in\mathbb{T}$, so 
$r_j\rightarrow 0\in\mathbb{T}$. Since $t_j\omega\rightarrow 0\in\mathbb{T}^d$ and $r_j\omega\rightarrow 0
\in\mathbb{T}^d$, 
then $n_j\omega\rightarrow 0\in\mathbb{T}^d$. By assumption, this implies that 
$n_j\beta\rightarrow 0\in 2\pi\mathbb{T}$. 
But $r_j\beta\rightarrow 0\in 2\pi\mathbb{T}$, so $t_j\beta\rightarrow 0\in 2\pi\mathbb{T}$. By 1., this implies that 
$\beta\in 2\pi\langle\mathbb{Z}^d,(\omega,1)\rangle$. $\Box$

\begin{lem}\label{alg} Let $W$ be a subspace of $\mathbb{C}^n$. Let $(W\cap \mathbb{R}^n)\otimes \mathbb{C}$ 
be the complex vector space generated by $W\cap \mathbb{R}^n$. Then 
\begin{enumerate}
\item $W=\bar{W}\Leftrightarrow (W\cap \mathbb{R}^n)\otimes \mathbb{C}=W$;
\item Let $V$ be a subspace of $W$ such that $V\oplus \bar{V}=W$ and $z_1,\dots z_k$ a basis of $V$, then 
$Re z_1,Im z_1,\dots Re z_k,Im z_k$ is a basis of $W\cap \mathbb{R}^n$.
\end{enumerate}
\end{lem}

\dem 
1. $\Rightarrow$: Let $w_1,\dots w_l$ be a basis of $W$; as $W=\bar{W}$, then $\bar{w_1},\dots \bar{w_l}$ 
are also in $W$. So, for all $j\in \{1,\dots l\}$, $Re w_j=\frac{1}{2}(w_j+\bar{w_j})$ and $Im w_j=\frac{1}{2i}
(w_j-\bar{w_j})$ are in $W\cap \mathbb{R}^n$. So $w_j=Re w_j+iIm w_j\in (W\cap \mathbb{R}^n)\otimes \mathbb{C}$ and
so $W\subset (W\cap \mathbb{R}^n)\otimes \mathbb{C}$. The other inclusion is obvious, therefore 
$W= (W\cap \mathbb{R}^n)\otimes \mathbb{C}$. 

\bigskip
$\Leftarrow$: Let $w\in (W\cap \mathbb{R}^n)\otimes \mathbb{C}$; then $w=\sum_{j=1}^la_jv_j$ where 
$a_j\in \mathbb{C}$ 
and $v_j\in  (W\cap \mathbb{R}^n)$. So $\bar{w}=\sum_{j=1}^l\bar{a_j}\bar{v_j}=\sum_{j=1}^l\bar{a_j}v_j\in 
(W\cap \mathbb{R}^n)\otimes \mathbb{C}$ and so $W=\bar{W}$. 

\bigskip
2. $Re z_1,Im z_1,\dots Re z_k,Im z_k$ generate $W$ and are real, so they generate $W\cap \mathbb{R}^n$. With complex
 coefficients, they generate $(W\cap \mathbb{R}^n)\otimes \mathbb{C}$, and by 1., this is equal to $W$. 
 As $W$ has dimension $2k$, they form a basis of $W\cap \mathbb{R}^n$. $\Box$

\subsection{Subbundles, invariant subbundles and Jordan subbundles}\label{sf}

\Defs 

\begin{itemize}

\item A real (resp. complex) subbundle is a family $V=\{V(\theta),\theta\in \mathbb{T}^d\}$ 
of subspaces of 
$\mathbb{R}^n$ (resp.$\mathbb{C}^n$) which is continuous in $\theta$, i.e such that for all 
$\theta_0\in\mathbb{T}^d$, there exists an open subset $\mathcal{U}$ containing $\theta_0$ and, for all 
$\theta\in\mathcal{U}$, 
a basis $\{x_1(\theta),\dots x_k(\theta)\}$ 
of $V(\theta)$ which is 
continuous in $\theta$ on $\mathcal{U}$. 

\item The dimension of $V(\theta)$ is automatically independant of $\theta$; 
this number is called the dimension of the subbundle $V$ and is denoted by $\mathrm{dim} V$.

\item A real (resp. complex) invariant subbundle for the cocycle $X$ is a real (resp. complex) subbundle 
such that for all $t,\theta$, $X^t(\theta)V(\theta)=V(\theta+t\omega)$. In what follows, we shall omit to mention 
the cocycle $X$, as no other cocycle is involved.

\end{itemize}

\bigskip
\rem A real invariant subbundle does not always have a basis which is continuous on $\mathbb{T}^d$.

\bigskip
\textbf{Example:} Consider the discrete 1-periodic cocycle $X^1(\theta):=\left(\begin{array}{cc}
\cos 2\pi\theta & -\sin 2\pi\theta\\
\sin 2\pi\theta & \cos 2\pi\theta\\
\end{array}\right)$ acting on $\mathbb{R}^2$. Let $z_1(\theta):=\left(\begin{array}{c}
\cos\pi\theta\\
\sin\pi\theta\\
\end{array}\right)$. Then $\mathrm{Vect}_\mathbb{R}(z_1(\theta))$ is an invariant subbundle for $X^n(\theta)$, 
since $z_1(\theta+1)=-z_1(\theta)$ for all $\theta$. But
 $z_1$ is continuous on $2\mathbb{T}$ and not on $\mathbb{T}$. Moreover, if $z$ is another function such that for all 
 $\theta$, $z(\theta)$ generates 
$\mathrm{Vect}_\mathbb{R}(z_1(\theta))$, then there exists a continuous function $\lambda$ 
bounded away from $0$ and such that for all $\theta$, $z(\theta)=\lambda(\theta)z_1(\theta)$. 
So $z(\theta+1)=\lambda(\theta+1)z_1(\theta+1)=-\lambda(\theta+1)z_1(\theta)$, and so $z(\theta)$ is continuous 
on $\mathbb{T}$ if and only if for all $\theta$, $-\lambda(\theta+1)=\lambda(\theta)$. But this implies that the 
function $\lambda$ changes sign, so it takes the value $0$ since it is continuous, which is impossible.

\bigskip
\rem The intersection of two (real or complex) subbundles is not necessarily a subbundle. For instance, in 
$\mathbb{R}^2$, 
for all $\theta\in \mathbb{T}$, let $V(\theta)=\left(\begin{array}{c}
1\\
0\\
\end{array}\right)$ et $W(\theta)=\left(\begin{array}{c}
\cos 2\pi \theta\\
\sin 2\pi \theta\\
\end{array}\right)$, then $V(\theta)\cap W(\theta)=V(\theta)$ if $\theta=0$ or $\frac{1}{2}\ \mathrm{mod}\ 1$ 
and $\{0\}$ otherwise, so the dimension of the intersection is not independant of $\theta$. 
However, the following proposition holds:

\begin{prop} The intersection of two real or complex invariant subbundles is an invariant subbundle. \end{prop}

\dem Let $U,V$ be two invariant subbundles, then for all $t,\theta$, 

$$X^t(\theta)(U(\theta)\cap V(\theta))
=X^t(\theta)U(\theta)\cap X^t(\theta)
V(\theta)=U(\theta+t\omega)\cap V(\theta+t\omega)$$ 

\noindent so the intersection is invariant.

\bigskip
Let us show that it has constant dimension. Let $\mathcal{U}$ be an open subset of the torus such that there exists 
$(u_1,
\dots u_k)$ and $(v_1,\dots v_l)$ continuous on $\mathcal{U}$ and such that for all $\theta\in\mathcal{U}$, 
$(u_1(\theta),\dots u_k(\theta))$ is a basis of $U(\theta)$ and $(v_1(\theta),\dots v_l(\theta))$ a basis of $V(\theta)$. 
 For all $\theta\in\mathcal{U}$, let
 
 $$M(\theta):=\left[\begin{array}{cccccc}
u_1(\theta) & \dots u_k(\theta) & v_1(\theta) & \dots & v_l(\theta)\\
\end{array}\right]$$ 

\noindent the $n\times (k+l)$-matrix whose columns are the vectors from the two bases. Let 
$r$ be the rank of $M(\theta_0)$ for a fixed $\theta_0$ in $\mathcal{U}$. 
Then there exists a $r\times r$-submatrix of 
$M(\theta_0)$ with non-zero determinant. This determinant is continuous in $\theta$, 
so it is non-zero on a neighbourhood 
$\mathcal{V}$ of $\theta_0$, so, on this neighbourhood, the rank of $M(\theta)$ is greater than, or equal to $r$. 
Therefore, if $d(\theta)$ is the dimension of $U(\theta)\cap 
V(\theta)$, then $d(\theta)\leq d(\theta_0)$ for all $\theta\in\mathcal{V}$. 

\bigskip
Let $\theta_0,\phi_0\in\mathbb{T}^d$. As we have just seen, there exists a neighbourood $\mathcal{U}_0$ of $\theta_0$ 
and a neighbourhood 
$\mathcal{V}_0$ of $\phi_0$ such that $d(\theta)\leq d(\theta_0)$ for all $\theta\in\mathcal{U}_0$ and 
$d(\phi)\leq d(\phi_0)$ for all 
$\phi\in\mathcal{V}_0$. As the orbits of $t\mapsto \theta_0+t\omega$ and $t\mapsto \phi_0+t\omega$ are dense on the 
torus, 
there exists $t,t'\in\mathbb{R}$ such that $\theta_0+t\omega\in\mathcal{V}_0$ and $\phi_0+t'\omega\in\mathcal{U}_0$. 
Invariance 
and invertibility of $X^t(\theta_0)$ imply that $d(\theta_0+t\omega)=\mathrm{dim} (U\cap V)(\theta_0+t\omega)
=\mathrm{dim} X^t(\theta_0) (U\cap V)(\theta_0)=\mathrm{dim}(U\cap V)(\theta_0)=d(\theta_0)$, and analogously 
$d(\phi_0)=d(\phi_0+t'\omega)$. 
Moreover, as $\theta_0+t\omega\in\mathcal{V}_0$, $d(\theta_0+t\omega)\leq d(\phi_0)$, and analogously 
$d(\phi_0+t'\omega)
\leq d(\theta_0)$. Therefore, $d(\theta_0)=d(\phi_0)$. As $\theta_0$ et $\phi_0$ are arbitrarily chosen, 
the dimension of 
$U\cap V$ is constant on $\mathbb{T}^d$. 

\bigskip
Let us now define a local basis of $U\cap V$. 
Let $\mathcal{U}$ be a sufficiently small neighbourhood of $\theta_0$ in $\mathbb{T}^d$ and 
$(u_1,\dots u_k)$ and 
$(v_1,\dots v_l)$ two bases for $U$ and $V$ which are continuous
on $\mathcal{U}$.
In the neighbourhood of $\theta_0$, up to a permutation of the bases, 
there exists $l'\leq l$ such that $u_1(\theta),\dots u_k(\theta), v_1(\theta),\dots v_{l'}(\theta)$ is a basis of 
$U(\theta)+V(\theta)$. 
the integer $l'$ does not depend on $\theta$ since the dimension of $U(\theta)\cap V(\theta)$ is independant of $\theta$. 
So, for all $l''$, $l'\leq l''\leq l$, there exists $a_1,\dots a_k,b_1,\dots b_{l'}$ which are continuous on a 
neignbourhood 
of $\theta_0$ such that $v_{l''}(\theta)=\sum_{i=1}^k a_i(\theta)u_i(\theta)+\sum_{i=1}^{l'}
b_i(\theta)v_i(\theta)$. 
Let $\bar{v}_{l''}(\theta):=\sum_{i=1}^k a_i(\theta)u_i(\theta)$, then $\bar{v}_{l''}(\theta)\in U(\theta)\cap V(\theta)$, 
$\bar{v}_{l''}$ is 
continuous on a neighbourhood of $\theta_0$. The vectors $(\bar{v}_{l'+1}(\theta_0),\dots \bar{v}_{l}(\theta_0))$ 
form a basis of 
$U(\theta_0)\cap V(\theta_0)$, therefore $\bar{v}_{l'+1}(\theta),\dots \bar{v}_l(\theta)$ form a basis of $U(\theta)
\cap V(\theta)$ in a neighbourhood of 
$\theta_0$. $\Box$

\bigskip
\Def A Jordan subbundle of rank $k$ modulo $N$ is a complex invariant subbundle having a basis $(z_1,\dots 
z_k)$ which is continuous on $N\mathbb{T}^d$ and such that there exists $\alpha+i\beta\in\mathbb{C}$ satisfying 
for all 
$\theta,t$,

\begin{equation}\begin{split}
& X^t(\theta)z_1(\theta)=e^{t(\alpha+i\beta)}z_1(\theta+t\omega)\\
& X^t(\theta)z_2(\theta)=e^{t(\alpha+i\beta)}z_2(\theta+t\omega)+te^{t(\alpha+i\beta)}z_1(\theta+t\omega)\\
& \dots\\
& X^t(\theta)z_k(\theta)=e^{t(\alpha+i\beta)}\sum_{i=1}^k\frac{t^{k-i}}{(k-i)!}z_i(\theta+t\omega)\\
\end{split}\end{equation}

\noindent A Jordan subbundle is a Jordan subbundle modulo 1. 
The family of functions 
$(z_1,
\dots z_k)$ is called a Jordan basis, it is not unique. If it is real for all $\theta$, it is called a real Jordan 
basis (for a complex Jordan subbundle). 
The number $\alpha+i\beta$ is called an exponent of the Jordan subbundle, and also the exponent of the Jordan basis 
$(z_1,
\dots z_k)$. 

\bigskip
\rem An exponent of a Jordan subbundle is not unique, but the exponent of a Jordan basis is. 

\noindent If unnecessary, we shall omit to mention the rank of a Jordan subbundle. Notice that the rank is not supposed to be 
maximal: if $k\geq 2$, a Jordan subbundle of rank $k$ contains another Jordan subbundle of rank $k-1$.

\bigskip
\Def An invariant subbundle $W$ of dimension $k$ is reducible modulo $N$ if there exists a basis 
$(z_1,\dots z_k)$ of $W$ which is continuous on $N\mathbb{T}^d$  
and a constant matrix $A$ of dimension $k\times k$ such that  
$X^t(\theta)[z_1(\theta)\dots z_k(\theta)]=[z_1(\theta+t\omega)\dots z_k(\theta+t\omega)]e^{tA}$ 
for all $t,\theta$.

\rem A Jordan subbundle is a particular type of reducible invariant subbundle and
 $GL(n,\mathbb{R})$-reducibility 
is equivalent to the existence of a decomposition of $\mathbb{R}^n$ into invariant reducible subbundles.

\begin{prop}\label{exponent} Let $V$ be a Jordan subbundle modulo $N$. 

i) If $\alpha+i\beta$ is an exponent for $V$, then for all $m\in\mathbb{Z}^d$, 
$\alpha+i\beta+2i\pi\langle\omega,\frac{m}{N}\rangle$ 
is an exponent for $V$.

ii) If $\alpha+i\beta$ and $\alpha'+i\beta'$ are two exponents for $V$, then $\alpha=\alpha'$ and 
$\beta-\beta'\in 2\pi\langle\mathbb{Z}^d,\frac{\omega}{N}\rangle$.
\end{prop} 

\dem i) Suppose $(z_1,\dots z_k)$ is a Jordan basis of $V$ with exponent $\alpha+i\beta$.

\noindent Let $m\in\mathbb{Z}^d$. 
For all $1\leq j\leq k$ and all $\theta\in N\mathbb{T}^d$, let $z'_j(\theta)=e^{-2i\pi\langle\frac{\theta}{N},m\rangle}z_j(\theta)$. Then the vectors 
$z'_j(\theta)$ form a global basis of $V$ which is 
continuous on $N\mathbb{T}^d$ and for all $\theta,t$ and all $j\leq k$,

\begin{equation}
X^t(\theta)z'_j(\theta)=e^{t(\alpha+i\beta+2i\pi\langle m,\frac{\omega}{N}\rangle)}\sum_{i=1}^j\frac{t^{j-i}}{(j-i)!}z'_i(\theta+t\omega)
\end{equation}

\noindent so $\alpha+i\beta+2i\pi\langle m,\frac{\omega}{N}\rangle$ is also an exponent of $V$.

\bigskip
ii) Let $(v_1,\dots v_k)$ and $(v'_1,\dots v'_k)$ be Jordan bases of $V$ with respective exponents 
$\alpha+i\beta$ and 
$\alpha'+i\beta'$.

\noindent For all $\theta\in N\mathbb{T}^d$, let $v'_1(\theta)=\sum_{j=1}^k\gamma_j(\theta)v_j(\theta)$ where $\gamma_j$ are continuous on 
$N\mathbb{T}^d$. 
Then for all $t$,

$$\sum_{j=1}^k\gamma_j(\theta)e^{t(\alpha+i\beta)}\sum_{i=1}^j\frac{t^{j-i}}{(j-i)!}v_i(\theta+t\omega)
=e^{t(\alpha'+i\beta')}\sum_{j=1}^k\gamma_j(\theta+t\omega)v_j(\theta+t\omega)$$

\noindent As the $v_j(\theta+t\omega)$ are linearly independant, in particular

$$\gamma_k(\theta)e^{t(\alpha+i\beta)}=e^{t(\alpha'+i\beta')}\gamma_k(\theta+t\omega)$$

\noindent Suppose $\gamma_k(\theta)\neq 0$ for some $\theta\in N\mathbb{T}^d$. 
As $\gamma_k$ is bounded, then $\alpha=\alpha'$.
Let $t_m$ be an unbounded real sequence such that $t_m\omega\rightarrow 0\in N\mathbb{T}^d$. 
Then as $m\rightarrow \infty$, 
since $\gamma_k(\theta)\neq 0$, $t_m({\beta}-{\beta'})\rightarrow 0\in 2\pi\mathbb{T}$. By lemma \ref{dep}, 
there exists $K\in\mathbb{Z}^d$ such that $\beta-\beta'=2\pi\langle K,\frac{\omega}{N}\rangle$.
If $\gamma_k$ is identically zero, then 

$$\gamma_{k-1}(\theta)e^{t(\alpha+i\beta)}=e^{t(\alpha'+i\beta')}\gamma_{k-1}(\theta+t\omega)$$

\noindent and we deduce in the same way that $\beta-\beta'=2\pi\langle K,\frac{\omega}{N}\rangle$ for some $K\in\mathbb{Z}^d$. 
Otherwise, we repeat the argument until we find a non zero $\gamma_j(\theta)$ and deduce that for some 
$K\in\mathbb{Z}^d$, 
$\beta-\beta'=2\pi\langle K,\frac{\omega}{N}\rangle$. $\Box$

\bigskip
\rem $\bullet$
Thus, the exponent of a Jordan subbundle modulo $N$ is well defined modulo $2i\pi\langle\mathbb{Z}^d,\frac{\omega}
{N}\rangle$. 
In particular, if $\beta\in 2\pi\langle\mathbb{Z}^d,\frac{\omega}{N}\rangle$, then we can assume that $\beta=0$.

\bigskip
$\bullet$ The term "Jordan subbundle" comes from the fact that if (\ref{ladef}) holds for some $B$ 
in Jordan normal form, then the columns of $Z(\theta)^{-1}$ whose indices are the same as those of the first columns 
of a Jordan block of $B$ with eigenvalue $\alpha+i\beta$ form a Jordan basis with exponent $\alpha+i\beta$.

\begin{lem}\label{rédsfJ} $GL(n,\mathbb{C})$-reducibility modulo $N$ is equivalent to the existence of a decmposition of 
$\mathbb{C}^n$ into Jordan subbundles modulo $N$. 
The existence of a decomposition of $\mathbb{R}^n$ into Jordan subbundles modulo $N$ with a real Jordan basis 
implies $GL(n,\mathbb{R})$-reducibility modulo $N$. \end{lem}

\dem By definition, $GL(n,\mathbb{C})$-reducibility of $X$ is the existence of a matrix 
$B=\left(\begin{array}{ccc}
B_1 & 0 & 0\\
0 & B_2 & 0\\
0 & 0 & \ddots\\
\end{array}\right)$, where each $B_j$ is a Jordan block with exponent $\alpha_j+i\beta_j$, and of a continuous function 
$Z:\mathbb{T}^d\rightarrow GL(n,\mathbb{C})$ such that for all $\theta,t$,

\begin{equation}\label{red2}X^t(\theta)=Z(\theta+t\omega)^{-1}e^{tB}Z(\theta)\end{equation}

\noindent If $z_1(\theta),\dots z_n(\theta)$ are the columns of $Z(\theta)^{-1}$, then (\ref{red2}) is equivalent to the fact that 
for all $\theta,t,j$, if $l_1,\dots l_{k_j}$ are the indices of the columns containing $B_j$, 

\begin{equation}
\begin{split}
& X^t(\theta)z_{l_1}(\theta)=e^{t(\alpha_j+i\beta_j)}z_{l_1}(\theta+t\omega)\\
& X^t(\theta)z_{l_2}(\theta)=e^{t(\alpha_j+i\beta_j)}z_{l_2}(\theta+t\omega)+te^{t(\alpha_j+i\beta_j)}z_{l_1}(\theta+t\omega)\\
& \dots\\
& X^t(\theta)z_{l_{k_j}}(\theta)=e^{t(\alpha_j+i\beta_j)}\sum_{i=1}^{k_j}\frac{t^{k_j-i}}{(k_j-i)!}z_i(\theta+t\omega)\\
\end{split}\end{equation}

\noindent which is also equivalent to the fact that if for all $j$, $V_j(\theta)=\mathrm{Vect}_\mathbb{C}(z_{l_1}(\theta),\dots z_{l_{k_j}}(\theta))$, 
then $V_j$ is a Jordan subbundle with exponent $\alpha_j+i\beta_j$. Moreover, $V_j(\theta)$ are in direct sum since 
$Z(\theta)^{-1}$ is invertible. 

\bigskip
In the preceding argument, it is clear that $X$ is in fact $GL(n,\mathbb{R})$-reducible if all 
$V_j$ have a real global basis.  
$\Box$

\bigskip
\rem $\bullet$ Decomposition into Jordan subbundles is not always unique. For instance, if for all $\theta,t$,  
$X^t(\theta)=Z(\theta+t\omega)^{-1}e^{t\alpha Id}Z(\theta)$, then for any invertible matrix $P$, 
$X^t(\theta)=Z(\theta+t\omega)^{-1}Pe^{t\alpha Id}P^{-1}Z(\theta)$, so $\mathbb{C}^n$ decomposes into Jordan subbundles 
of rank 1 generated by the columns of $Z(\theta)^{-1}P$, where $P$ is arbitrarily chosen.

\bigskip
$\bullet$ If a Jordan subbundle has a real Jordan basis with exponent ${\alpha+i\beta\ 
(\mathrm{mod} 2i\pi\langle\mathbb{Z}^d,
\omega\rangle)}$, then $\beta=0\ (\mathrm{mod} 2\pi\langle\mathbb{Z}^d,
\omega\rangle)$. 
But there exists Jordan subbundles with real exponent $(\mathrm{mod} 2i\pi\langle\mathbb{Z}^d,
\omega\rangle)$ but without a real Jordan basis. A trivial example 
is the constant Jordan subbundle generated by $\left(\begin{array}{c}
1\\
i\\
\end{array}\right)$ with exponent $0\ (\mathrm{mod} 2i\pi\langle\mathbb{Z}^d,
\omega\rangle)$ for the identity cocycle. However, the following lemma holds:

\begin{lem}\label{baseU} Let $W$ be a Jordan subbundle modulo N without a real Jordan basis, 
$z_1,\dots z_k$ a Jordan basis of $W$ with real exponent $\alpha$, and for all $j\leq k$, $u_j$ the real part of $z_j$ and $v_j$ its imaginary part. 
Then $U:=\mathrm{Vect}_\mathbb{R}(u_1,\dots u_k)$ and $V:=\mathrm{Vect}_\mathbb{R}(v_1,\dots v_k)$ are Jordan 
subbundles modulo N with exponent $\alpha$ and there exists $l,m\leq k$ such that $u_l,\dots u_k$ is a Jordan basis 
of $U$ and 
$v_m,\dots v_k$ a Jordan basis of $V$. Moreover, either $l$ or $m$ is equal to 1.\end{lem}

\dem For all $j$, $t\in\mathbb{R},\theta\in N\mathbb{T}^d$, as $X^t(\theta)$ is real, then 

$$X^t(\theta)u_j(\theta)=e^{t\alpha}\sum_{j'\leq j}\frac{t^{j-j'}}{(j-j')!}u_{j'}(\theta+t\omega)$$

\noindent Suppose there exists $j\geq 1,\theta_0$ and $\lambda_1,\dots \lambda_{j-1}\in\mathbb{C}$ such that 
$u_j(\theta_0)=\sum_{i\leq j-1}\lambda_i u_i(\theta_0)$. 
Then for all $t$,

\begin{equation}\label{rgjo}\begin{split}&0=X^t(\theta_0)(u_j(\theta_0)-\sum_{i\leq j-1}\lambda_i u_i(\theta_0))\\
&=e^{t\alpha}\sum_{j'\leq j}\frac{t^{j-j'}}{(j-j')!}u_{j'}(\theta_0+t\omega)
-\sum_{i\leq j-1}\lambda_i \sum_{j'\leq i}\frac{t^{i-j'}}
{(i-j')!} u_{j'}(\theta_0+t\omega)
\end{split}
\end{equation}

\noindent so, dividing by $e^{t\alpha}t^{j-1}$, for all $t\neq 0$,

$$0=\sum_{j'\leq j}\frac{t^{-j'+1}}{(j-j')!}u_{j'}(\theta_0+t\omega)-\sum_{i\leq j-1}\lambda_i 
\sum_{j'\leq i}\frac{t^{i-j'-j+1}}
{(i-j')!} u_{j'}(\theta_0+t\omega)$$

\noindent Let $\theta$ be any point of $N\mathbb{T}^d$. 
Let $t_s$ be an unbounded real sequence satisfying $t_s\omega\rightarrow \theta-\theta_0$ in $N\mathbb{T}^d$. 
Then, as $s$ tends to infinity,

$$\frac{1}{(j-1)!}u_1(\theta)=0$$

\bigskip
Assume by induction that $u_{j''}$ is identically $0$ for all $j''$ 
strictly inferior to some $J\leq j$. Then, dividing equation (\ref{rgjo}) by $e^{t\alpha}t^{j-J}$, 
for all $t\neq 0$, 

\begin{equation}0=e^{t\alpha}\sum_{J\leq j'\leq j}\frac{t^{J-j'}}{(j-j')!}u_{j'}(\theta_0+t\omega)
-\sum_{i\leq j-1}\lambda_i \sum_{J\leq j'\leq i}\frac{t^{J-j+i-j'}}
{(i-j')!} u_{j'}(\theta_0+t\omega))\end{equation}

\noindent so, with the sequence $t_s$ above defined, if $t=t_s$ and taking the limit as $s\rightarrow \infty$, 

$$\frac{1}{(j-J)!}u_J(\theta)=0$$

\noindent and so $u_J(\theta)=0$ for all $\theta$. Therefore, for all $\theta$ and all $j'\leq j$, $u_{j'}(\theta)=0$.

\bigskip
Thus, we have shown that there exists $l\leq k$ so that the functions $u_1,\dots u_{l-1}$ are identically $0$ 
if $l\geq 2$ 
and $(u_l,\dots u_k)$ form a global basis of $U$, which is then a Jordan basis. 
We proceed exactly in the same way to show that there exists 
$m\leq k$ such that $v_1,\dots v_{m-1}$ are identically $0$ if $m\geq 2$ 
and $(v_m,\dots v_k)$ form a Jordan basis of $V$. Moreover, as $u_1$ and $v_1$ cannot be $0$ at the same time, 
then either $l$ or $m$ is equal to 1. $\Box$

\subsection{Properties of Jordan subbundles with a real Jordan basis}

Let $\{ u_j,1\leq j\leq k\}$ be a real Jordan basis  
of a Jordan subbundle $U$ modulo N of rank $k$ and real exponent $\alpha$. 

\begin{soulem}\label{blabla} Every invariant subbundle contained in $U$ is a Jordan subbundle modulo N 
generated by 
$(u_1,\dots u_j)$ for some 
$j\leq k$.\end{soulem}

\dem 
Let $W$ be a non zero invariant subbundle contained in $U$ and $u_1,\dots u_k$ as above. 
For some $\theta_0$, let $j$ be the maximal integer lower than $k$ such that there exists 
$\sum_{j'\leq j} a_{j'}u_{j'}(\theta_0)$ in $W(\theta_0)$ with $a_j\neq 0$. 

\noindent As $W$ is invariant, for all $t\in\mathbb{R}$, $X^t(\theta_0)\sum_{j'\leq j}a_{j'}
u_{j'}(\theta_0)\in W(\theta_0+t\omega)$. 
Now this vector is equal to 
$e^{t\alpha}\sum_{j'\leq j}a_{j'}\sum_{i=1}^{j'}\frac{t^{j'-i}}{(j'-i)!}u_i(\theta_0+t\omega)$.
Dividing by $e^{t\alpha}t^{j-1}$, for all $t\neq 0$, the vector $\sum_{j'\leq j}a_{j'}\sum_{i=1}^{j'}\frac{t^{j'-j+1-i}}
{(j'-i)!}u_{i}(\theta_0
+t\omega)$ is in $W(\theta_0+t\omega)$. 
Let $\theta\in N\mathbb{T}^d$ and $(t_k)$ an unbounded real sequence such that $t_k\omega\rightarrow \theta
-\theta_0$ in $N\mathbb{T}^d$. 
Then, taking the limit as $k\rightarrow \infty$, 

$$\sum_{j'\leq j}a_{j'}\sum_{i=1}^{j'}\frac{t_k^{j'-j+1-i}}{(j'-i)!}u_{i}(\theta_0+t_k\omega)\rightarrow 
\frac{1}{(j-1)!}a_ju_1(\theta)
\in W(\theta)$$

\noindent So for all $\theta$, $u_1(\theta)\in W(\theta)$. Suppose that for all $1\leq j''\leq j'$, $u_{j''}(\theta)
\in W(\theta)$ for all $\theta$. Then $\sum_{i=j'+1}^{j}a_iu_i(\theta_0)\in W(\theta_0)$ so by invariance of $W$, 
for all $k$,

$$\sum_{i=j'+1}^{j}e^{t_k\alpha}\sum_{j''=1}^ia_i\frac{t_k^{i-j''}}{(i-j'')!}u_{j''}(\theta_0+t_k\omega)
\in W(\theta_0+t_k\omega)$$

\noindent and so

$$\sum_{i=j'+1}^{j}e^{t_k\alpha}\sum_{j''=j'+1}^ia_i\frac{t_k^{i-j''}}{(i-j'')!}u_{j''}(\theta_0+t_k\omega)
\in W(\theta_0+t_k\omega)$$

\noindent Dividing by $e^{t_k\alpha}t_k^{j-j'-1}$, we get

$$\sum_{i= j'+1}^ja_{i}\sum_{j''=j'+1}^{i}\frac{t_k^{i-j''+j'-j+1}}
{(i-j'')!}u_{j''}(\theta_0+t_k\omega)\in W(\theta_0+t_k\omega)$$

\noindent and taking the limit as $k$ goes to infinity, as $a_j\neq 0$, $u_{j'+1}(\theta)\in W(\theta)$. 
Eventually, $\mathrm{Vect}_\mathbb{C}(u_1,\dots u_j)$ 
is contained in $W$; therefore, since we assumed $j$ is maximal, $\mathrm{Vect}_\mathbb{C}(u_1,\dots u_j)$ 
is equal to $W$. $\Box$

\begin{soulem}\label{bla} Let $W'$ be an invariant subbundle such that for all $\theta\in \mathbb{T}^d$, $W'(\theta)=
\bigoplus_{i=1}^mW^i(\theta)$ where each $W^i$ is a Jordan subbundle modulo N with a real 
Jordan basis 
$(w^i_1,\dots w^i_{l_i})$ with exponent $\alpha$ and suppose $u_j(\theta_0)=\sum_{i=1}^m\sum_{l=1}^{l_i}\lambda^l_i w^i_l(\theta_0)$ 
for some $\theta_0$, then for all $\theta\in N\mathbb{T}^d$ and all $j'\leq j$,  
$u_{j'}(\theta)=\sum_{i=1}^m\sum_{1\leq j''\leq \min(j',l_i-j+j')}\lambda_i^{j-j'+j''} w^i_{j''}(\theta)$.
 In particular, $u_{1}(\theta)=\sum_{i=1}^m\lambda_i^{j} w^i_{1}(\theta)$.
\end{soulem}

\dem For all $t$,

\begin{equation}\label{gf}\begin{split}&0=X^t(\theta_0)(u_j(\theta_0)-\sum_{i=1}^m\sum_{l=1}^{l_{i}}\lambda^l_{i}
 w^{i}_l(\theta_0))\\
&=e^{t\alpha}
\sum_{j'\leq j}\frac{t^{j-j'}}{(j-j')!}u_{j'}(\theta_0+t\omega)-\sum_{i=1}^m\sum_{l=1}^{l_i}\lambda^l_{i}
e^{t\alpha}
\sum_{l'=1}^l\frac{t^{l-l'}}{(l-l')!} w^{i}_{l'}(\theta_0+t\omega)\end{split}\end{equation}

\noindent Dividing by $e^{t\alpha}$, we get for all $t$

\begin{equation}\label{grrrr}0=
\sum_{i\leq j}\frac{t^{j-j'}}{(j-j')!}u_{j'}(\theta_0+t\omega)-\sum_{i=1}^m\sum_{l=1}^{l_i}\lambda^l_{i} 
\sum_{l'=1}^l\frac{t^{l-l'}}{(l-l')!}w^{i}_{l'}(\theta_0+t\omega)\end{equation}

\noindent Let $L$ be the greatest power of $t$ in this expression. 
Let $\theta$ be any point of $N\mathbb{T}^d$. 
Take a sequence $t_k\rightarrow \infty$ such that $t_k\omega\rightarrow \theta-\theta_0\in N\mathbb{T}^d$ as 
$k\rightarrow \infty$. Suppose first that 
$L\geq j$. Then, dividing (\ref{grrrr}) by $t^L$, and making $k$ go to infinity,

\begin{equation}\sum_{i=1}^m\sum_{l=L+1}^{l_i}\lambda^l_{i} \frac{1}{L!}
w^{i}_{l-L}(\theta)= 0\end{equation}

\noindent Since $w^i_{l-L}(\theta)$ are linearly independant, $\lambda^l_i=0$ if $l\geq L+1$. Consequently, 
(\ref{grrrr}) can be rewritten

\begin{equation}0=\sum_{j'\leq j}\frac{t^{j-j'}}{(j-j')!}u_{j'}(\theta_0+t\omega)
-\sum_{i=1}^m \sum_{l=1}^{L}\lambda^l_{i} 
\sum_{l'=1}^l\frac{t^{l-l'}}{(l-l')!}w^{i}_{l'}(\theta_0+t\omega)
\end{equation}

\noindent But this contradicts the definition of $L$, so the assumption under which $L\geq j$ is false. Therefore, 
(\ref{grrrr}) can be rewritten

\begin{equation}\label{grrr}0=\sum_{j'\leq j}\frac{t^{j-j'}}{(j-j')!}u_{j'}(\theta_0+t\omega)
-\sum_{i=1}^m\sum_{l=1}^{\min(j,l_i)}\lambda^l_{i}\sum_{l'=1}^l \frac{t^{l-l'}}{(l-l')!}w^{i}_{l'}(\theta_0+t\omega)
\end{equation}

\noindent Dividing (\ref{grrr}) 
by $t^{j-1},\dots t$, replacing $t$ by $t_k$ and making $k$ go to $\infty$, we see that for all $1\leq j'\leq j$ 
and all $\theta\in N\mathbb{T}^d$,
 
\begin{equation}\begin{split}&u_{j'}(\theta)=\sum_{i=1}^m\sum_{l-l'=j-j'}\lambda^l_{i}w^{i}_{l'}(\theta)
=\sum_{i=1}^m\sum_{l=j-j'+1}^{\min(j,l_i)}\lambda^{l}_{i}w^{i}_{l-j+j'}(\theta)\\
&=\sum_{i=1}^m\sum_{l=1}^{\min(j',l_i-j+j')}\lambda^{l+j-j'}_{i}w^{i}_{l}(\theta)\ \Box\end{split}\end{equation}

\bigskip
\rem Coefficients 
$\lambda^{j-j'+j''}_i$ do not depend on $\theta$.

\begin{lem}\label{blablabla}
Let $W'$ be an invariant subbundle such that for all $\theta\in \mathbb{T}^d$, $W'(\theta)=\bigoplus_{i'=1}^mW^{i'}(\theta)$ 
where $W^{i'}$ are Jordan subbundles modulo N with a real basis $(w_1^{i'},
\dots w^{i'}_{l_{i'}})$ with exponent $\alpha$. 
Then $W'+U$ is a direct sum of Jordan subbundles modulo N with a real basis.\end{lem}

\dem If $U(\theta)\cap W'(\theta)=\{0\}$ for all $\theta$, this is trivial.

\noindent Let us now suppose that this intersection is non trivial. 
It is then equal to some non trivial invariant subbundle. 
By sublemma \ref{blabla}, it is generated by
 $u_1,\dots u_j$ for some $j\leq k$. 

\bigskip
Assume first that $\mathrm{dim}U\leq \mathrm{dim}W_i$ for all $i$. 

\noindent By sublemma \ref{bla}, there exists $\lambda_1, \dots \lambda_m$ such that for all $\theta\in N\mathbb{T}^d$,

$$u_1(\theta)=\sum_{i=1}^m\lambda_i w_1^i(\theta)$$

\noindent Let $u'_1=u_2-\sum_{i=1}^m\lambda_i w_2^i$,...$u'_{n-1}=u_n-\sum_{i=1}^m\lambda_i w_n^i$. 
If $u'_1,\dots u'_{n-1}$ are a basis, since for all $j\leq n-1$

\begin{equation}\begin{split}&
X^t(\theta)u'_j(\theta)=X^t(\theta)(u_{j+1}(\theta)-\sum_{i=1}^m\lambda_i w_{j+1}^i(\theta))\\
&=\sum_{j'\leq j+1}\frac{t^{j-j'}}{(j-j')!}(u_{j'}(\theta+t\omega)-\sum_{i=1}^m\lambda_i w_{j'}^i(\theta+t\omega))\\
& =\sum_{j'\leq j} \frac{t^{j-j'}}{(j-j')!}u'_{j'}(\theta+t\omega)
\end{split}
\end{equation}

\noindent this means that they are a Jordan basis.

\noindent If $u_1'$ is in the space generated by $u_2',\dots u_{n-1}'$, then we carry out the same construction. 
After finitely many steps, we have defined a Jordan basis for $U+\bigoplus_i W_i$.

\bigskip
Let now $U$ be of any dimension. We shall proceed by induction. 

\begin{itemize}
\item If $U$ has dimension 1, it is included in $W'$ so the conclusion immediatly follows.

\item Suppose now that the conclusion holds for any $U$ of dimension $\leq n-1$. 
If now $U$ has dimension $n$, write
$W'= W_1\oplus W_2$ où $W_1=\bigoplus_{dim W_i<n}W_i$ and $W_2=\bigoplus_{dim W_i\geq n}W_i$. By the above, 
$W_2+U$ is the direct sum of Jordan subbundles modulo $N$ with a real basis. Then, we add one by one the 
$W_i$ with dimension $<n$, and by induction hypothesis we still get a direct sum of Jordan subbundles modulo $N$ 
with a real basis. $\Box$

\end{itemize}

\subsection{Decomposition into invariant subbundles}\label{decomp}

Suppose $X$ is $GL(n,\mathbb{C})$-reducible. Then by lemma \ref{rédsfJ}, 

$$\forall \theta\in \mathbb{T}^d,\ 
\mathbb{C}^n=W_1(\theta)\oplus \dots \oplus W_r(\theta)$$ 

\noindent where each $W_j$ is the sum of all the Jordan subbundles 
with exponent ${\alpha_j+i\beta_j\ (\mathrm{mod} 2i\pi\langle\mathbb{Z}^d,
\omega\rangle)}$.

\begin{lem}\label{wj} For all $1\leq j\leq r$, there exists $1\leq j'\leq r$ such that $\bar{W_j}=W_{j'}$. 
Moreover, $W_j=\bar{W}_j$ iff $\beta_j\in \pi\langle\mathbb{Z}^d,\omega\rangle$.
\end{lem}

\dem Let $v(\theta)\in W_j(\theta)$ generating a complex invariant subbundle of dimension 1, 
then for all $\theta\in \mathbb{T}^d,t\in\mathbb{R}$,

$$X^t(\theta)v(\theta)=e^{t(\alpha_j+i\beta_j)}v(\theta+t\omega)$$ 

\noindent and

$$X^t(\theta)\bar{v}(\theta)=e^{t(\alpha_j-i\beta_j)}\bar{v}(\theta+t\omega)$$ 

\noindent Write $\bar{v}(\theta)=
\sum_{l=1}^r\gamma_l(\theta){w_l}(\theta)$ with $(w_l(\theta))_{l=1\dots r}$ a Jordan basis of $\mathbb{C}^n$ 
and $\gamma_l$ continuous and $\mathbb{C}$-valued. 
Then there exist polynomials $\{P_l(t),l=1,\dots r\}$ such that for all $t$,

$$X^t(\theta)\bar{v}(\theta)=\sum_{l=1}^r\gamma_l(\theta)e^{t(\alpha_l+i\beta_l)}P_l(t){w_l}(\theta+t\omega)$$ 

\noindent So

$$e^{(\alpha_j-i\beta_j)t}\sum_{l=1}^r\gamma_l(\theta){w_l}(\theta+t\omega)
=\sum_{l=1}^r\gamma_l(\theta)e^{(\alpha_l+i\beta_l)t}P_l(t){w_l}(\theta+t\omega)$$
 
\noindent Since $w_l(\theta+t\omega)$ are linearly independant, for all $l$,

$$e^{(\alpha_j-i\beta_j)t}\gamma_l(\theta)w_l(\theta+t\omega)
=\gamma_l(\theta)e^{(\alpha_l+i\beta_l)t}P_l(t)
w_l(\theta+t\omega)$$ 

\noindent so if $\gamma_l(\theta)\neq 0$, then for all $t$

$$w_l(\theta+t\omega)=e^{(\alpha_l-\alpha_j+i(\beta_j+\beta_l))t}P_l(t)w_l(\theta+t\omega)$$ 

\noindent This implies that $\alpha_l=\alpha_j$, $P_l$ is constant equal to 1 and $\beta_j=-\beta_l$. 

\noindent Let $j'$ be such that $w_l(\theta)\in W_{j'}(\theta)$ for all $\theta\in \mathbb{T}^d$; then $\bar{W_j}=W_{j'}$. 

\bigskip
Suppose now that $W_j=\bar{W}_j$. Let $V_1,\dots V_{R_j}$ be the Jordan subbundles contained in $W_j$, 
and for each $V_s$, $u_1^s+iv_1^s,\dots 
u^s_{k_s}+iv^s_{k_s}$ a global basis with exponent $\alpha+i\beta$. 
Write for all $\theta$ the decomposition $u^s_1(\theta)-iv^s_1(\theta)
=\sum_{s'\leq r,j\leq k_{s'}} a^{s'}_j(\theta)(u^{s'}_j(\theta)
+iv^{s'}_j(\theta))$, then let $X^t(\theta)$ act on each side; then for all $t$,

\begin{equation}\begin{split} & X^t(\theta)(u^s_1(\theta)-iv^s_1(\theta))=e^{t(\alpha-i\beta)}
(u^s_1(\theta+t\omega)-iv^s_1(\theta+t\omega))\\
&= e^{t(\alpha-i\beta)} \sum_{s'\leq r,j\leq k_{s'}} a^{s'}_j(\theta+t\omega)(u^{s'}_j(\theta+t\omega)
+iv^{s'}_j(\theta+t\omega)) \\
&= \sum_{s'\leq r,j\leq k_{s'}} a^{s'}_j(\theta)X^t(\theta)(u^{s'}_j(\theta)
+iv^{s'}_j(\theta))\\
&= \sum_{s'\leq r,j\leq k_{s'}} a^{s'}_j(\theta) e^{t(\alpha+i\beta)}\sum_{j'\leq j} \frac{t^{j-j'}}
{(j-j')!} (u^{s'}_{j'}(\theta+t\omega)
+iv^{s'}_{j'}(\theta+t\omega))\end{split}\end{equation}

\noindent as $u_{k_s}^s(\theta+t\omega)+iv_{k_s}^s(\theta+t\omega)$ is linearly independant from the rest, then 

$$e^{t(\alpha-i\beta)} a^{s}_{k_s}(\theta+t\omega)= a^{s}_{k_s}(\theta)e^{t(\alpha+i\beta)}$$

\noindent whence, by lemma \ref{dep}, the fact that $2\beta= 2\pi\langle m,\omega\rangle$ for some $m\in\mathbb{Z}^d$.

\bigskip
Conversely, if $2\beta= 2\pi\langle m,\omega\rangle$ for some $m\in\mathbb{Z}^d$, 
then $W_j$ is it own complex conjugate.
$\Box$

\subsection{Main result}\label{pp}

We get to the proof of theorem \ref{CR'}.

\begin{prop}\label{CR} Assume that the continuous cocycle $X$ is 
$GL(n,\mathbb{C})$-reducible. Then there exists a decomposition of $\mathbb{R}^n$ into two invariant subbundles 
$\mathcal{W}$ and $\mathcal{W}'$ 
such that:
\begin{itemize}
\item $\mathcal{W}$ is a reducible subbundle modulo 2, generated by a basis $(z_1,\dots z_s)$ such that for all 
$(\theta,t)\in 2\mathbb{T}^d\times\mathbb{R}$, 
$X^t(\theta)[z_1(\theta)\dots z_s(\theta)]=[z_1(\theta+t\omega)\dots z_s(\theta+t\omega)]e^{A_1t}$ where 
$A_1$ has a real spectrum;
\item $\mathcal{W}'$ is a reducible subbundle modulo 1 with a basis $(z_{s+1},\dots z_n)$ such that for all $(\theta,t)
\in\mathbb{T}^d\times \mathbb{R}$,

$${X^t(\theta)[z_{s+1}(\theta)\dots z_n(\theta)]=[z_{s+1}(\theta+t\omega)\dots z_n(\theta+t\omega)]e^{A_2t}}$$
 
\noindent with $\sigma(A_2)\cap (\mathbb{R}+i\pi\langle\mathbb{Z}^d,\omega\rangle\setminus\{0\})=\emptyset$ and if $\alpha_1+i\beta_1,\alpha_2+i\beta_2
\in \sigma(A_2)$, then $\beta_1-\beta_2$ is not in $2\pi\langle\mathbb{Z}^d,\omega\rangle\setminus\{0\}$.
\end{itemize}

\end{prop}

\dem
From lemma \ref{rédsfJ}, we get a decomposition of $\mathbb{C}^n$ 
into complex Jordan subbundles.

\noindent Let us keep the notations introduced in section \ref{decomp}.

\bigskip
By lemma \ref{wj}, there exists a decomposition $\mathbb{C}^n=W\oplus W'$ where $W$ is the direct sum of all $W_j$ 
which are their own complex conjugate, $W=\bigoplus_{j=1}^{r'}W_j$, and $W'$ the direct sum of all the others: 
$W'=\bigoplus_{j=r'+1}^rW_j$.  

\bigskip 
1. $\bullet$ By lemma \ref{wj}, 
$W$ contains exactly the Jordan subbundles whose exponent is in $\mathbb{R}+i\pi\langle\mathbb{Z}^d,\omega\rangle$. 
Decompose again $W$ into $W_\mathbb{R}$ and $W_\mathbb{C}$ where $W_\mathbb{R}$ is the sum of the Jordan 
subbundles having exponent $0\ mod\ 2i\pi\langle \mathbb{Z}^d,\omega\rangle$, and $W_\mathbb{C}$ is the sum 
of the Jordan 
subbundles whose exponent is in $i\pi\langle \mathbb{Z}^d,\omega\rangle\setminus 2i\pi\langle \mathbb{Z}^d,\omega\rangle$. 

\bigskip
For $W_\mathbb{R}$, we can find real Jordan bases with real exponent which are continuous on $\mathbb{T}^d$. 
 
\bigskip
Proposition \ref{exponent} implies that we can find real exponents for $W_\mathbb{C}$, but for bases which are 
continuous on $2\mathbb{T}^d$ and not on $\mathbb{T}^d$ anymore.

\bigskip
$\bullet$ We will show by induction that there is a decomposition of each $W_j\subset W_\mathbb{C}$ 
into Jordan subbundles 
with a real Jordan basis. 

\noindent Let $V_1,\dots V_{R_j}$ be the Jordan subbundles included in $W_j$. 
According to lemmas \ref{baseU} and \ref{blablabla}, 
$(V_1+ \bar{V_1})\cap \mathbb{R}^n$ is the direct sum of two Jordan subbundles modulo 2 with a real basis, 
since it is the sum of the Jordan subbundle modulo 2 generated by the real parts of the vectors in the basis of 
$V_1$,  
and of the Jordan subbundle modulo 2 generated by their imaginary parts.
 
\bigskip
Let $\bar{W}$ and $\bar{W}'$ be invariant subbundles such that there exists $k\geq 2$ with
$\bar{W}=(V_k+\bar{V}_k)\cap\mathbb{R}^n$ and that $\bar{W}'$ is a direct sum of Jordan subbundles modulo 2 
with a real basis. 

\noindent By lemma \ref{blablabla}, $\bar{W}$ is the direct sum of two Jordan subbundles modulo 2, $U$ and $V$. 

\noindent Using lemma \ref{blablabla} again, 
$U+\bar{W}'$ is the direct sum of Jordan subbundles modulo 2 with a real Jordan basis. 

\noindent Finally, by lemma \ref{blablabla}, $\bar{W}+\bar{W}'=V+U+\bar{W}'$ is the direct sum of Jordan subbundles 
modulo 2 with a real basis, 
which ends the induction.

\bigskip
2. In $W'$, choose for all $r'+1\leq j\leq r$ and for each Jordan subbundle  
$V^j_s,s\leq R_j$ contained in some $W_j\subset W'$, 
a Jordan basis with exponent $\alpha_j+i\beta_j$ 
such that for all $j,j'$, $\beta_j-\beta_{j'}$ is not in $2\pi\langle\mathbb{Z}^d,\omega\rangle\setminus\{0\}$. 
We have already showed that for all $j$, $\beta_j$ is not in $\pi\langle\mathbb{Z}^d,\omega\rangle$. 

\noindent Let $W''$ be a sum of Jordan subbundles such that $W'=W''\oplus \bar{W}''$. 
If $(u_1+iv_1,\dots u_{\frac{S}{2}}+iv_{\frac{S}{2}})$ is the global basis of $W''$ which is the union of all those 
Jordan bases, then lemma \ref{alg} implies that $(u_{1}(\theta),v_{1}(\theta),
\dots u_{\frac{S}{2}}(\theta),
v_{\frac{S}{2}}(\theta))$ form a basis of $W'(\theta)\cap \mathbb{R}^n$ for all $\theta$. Moreover, 

$$X^t(\theta)\big[u_{1}(\theta)\ v_{1}(\theta)\ 
\dots u_{\frac{S}{2}}(\theta)\ 
v_{\frac{S}{2}}(\theta)\big]=[u_{1}(\theta+t\omega)\ v_{1}(\theta+t\omega)\ 
\dots u_{\frac{S}{2}}(\theta+t\omega)\ 
v_{\frac{S}{2}}(\theta+t\omega)\big]e^{tA_2}$$ 

\noindent where $\sigma(A_2)=\{\alpha_j+i\beta_j,r'+1\leq j\leq r\}$. 

\bigskip
Let $\mathcal{W}=W_\mathbb{C}\cap \mathbb{R}^n$ and $\mathcal{W}'=W_\mathbb{R}\cap\mathbb{R}^n\oplus W'\cap 
\mathbb{R}^n$. We have shown the existence of the required bases $(z_1,\dots z_s)$ for $\mathcal{W}$ and 
$(z_{s+1},\dots z_n)$ for $\mathcal{W}'$.
$\Box$

\bigskip
\begin{cor}\label{youpi} With the notations of the proposition \ref{CR}, let 
$Z(\theta)=\big(z_1(\theta)\dots z_n(\theta)\big)$. Then for all $\theta,t$,

$$X^t(\theta)=Z(\theta+t\omega)e^{t\left(\begin{array}{cc}
A_1 & 0\\
0 & A_2\\
\end{array}\right)}Z(\theta)^{-1}$$ 

\end{cor}

\bigskip
This proves theorem \ref{CR'}.

\section{Reducibility in other Lie groups}

We now give the proof of the reducibility theorem for the groups $SL(n,\mathbb{R})$, 
$Sp(n,\mathbb{R})$, $O(n)$ and $U(n)$.

\subsection{$SL(n,\mathbb{R})$-reducibility}

\begin{prop}\label{RSL}

Let $X$ be a continuous $SL(n,\mathbb{R})$-valued cocycle which is $GL(n,\mathbb{R})$-reducible modulo $N$ 
to a cocycle $t\mapsto e^{tB}$.
Then $B\in sl(n,\mathbb{R})$ and there exists $\tilde{Z}:N\mathbb{T}^d\rightarrow SL(n,\mathbb{R})$ such that for all 
$t,\theta$,

\begin{equation}X^t(\theta)=\tilde{Z}(\theta+t\omega)^{-1}e^{tB}\tilde{Z}(\theta)\end{equation}

\noindent so $X^t(\theta)$ is $SL(n,\mathbb{R})$-reducible modulo $N$. 
\end{prop}

\dem Let $\tilde{Z}(\theta):=\frac{1}{\mathrm{det}Z(\theta)}Z(\theta)$. By construction, $\tilde{Z}\in C^0 
(N\mathbb{T}^d,SL(n,\mathbb{R}))$ and for all $\theta,t$,

\begin{equation}e^{tB}=Z(\theta+t\omega)X^t(\theta)Z(\theta)^{-1}\end{equation}

\noindent so

\begin{equation}\frac{\mathrm{det}Z(\theta)}{\mathrm{det}Z(\theta+t\omega)}e^{tB}=
\tilde{Z}(\theta+t\omega)X^t(\theta)\tilde{Z}(\theta)^{-1}\end{equation}

\noindent Thus the left-hand side has determinant 1. So

\begin{equation}\forall  t,\ \mathrm{tr}(\ln \frac{\mathrm{det}Z(\theta)}{\mathrm{det}Z(\theta+t\omega)}I+tB)=0\end{equation}

\noindent As $\ln \frac{\mathrm{det}Z(\theta)}{\mathrm{det}Z(\theta+t\omega)}$ is bounded, then $\mathrm{tr}(B)=0$ and 
$\mathrm{det}Z$ is constant. Therefore, for all $\theta,t$,

\begin{equation}X^t(\theta)=Z(\theta+t\omega)^{-1}\mathrm{det}Ze^{tB}\frac{Z(\theta)}{\mathrm{det}Z}
=\tilde{Z}(\theta+t\omega)^{-1}e^{tB}\tilde{Z}(\theta)\end{equation}

\begin{cor} Let $X$ be a $GL(n,\mathbb{C})$-reducible cocycle. 
If it is $SL(n,\mathbb{R})$-valued then it is $SL(n,\mathbb{R})$-reducible modulo 2.

\end{cor}

\dem Apply proposition \ref{CR}, then proposition \ref{RSL}. $\Box$

\bigskip
This proves theorem \ref{Sp'etc} when $G=SL(n,\mathbb{R})$.

\subsection{Symplectic reducibility}

\begin{prop}\label{CSp} If $X$ is $Sp(2n,\mathbb{R})$-valued and $GL(2n,\mathbb{C})$-reducible, 
then it is $Sp(2n,\mathbb{R})$-reducible modulo 2.\end{prop}

\dem 
Let $\mathbb{R}^n=\mathcal{W}\oplus \mathcal{W}'$ as in proposition \ref{CR} and $Z$ as in corollary  
\ref{youpi}: for all $\theta$, $Z(\theta)=[z_1(\theta)\ \dots z_n(\theta)]$.

\noindent Write $X^t(\theta)={Z}(\theta+t\omega)Ce^{tB}C^{-1}Z(\theta)^{-1}$ 
with $B$ in Jordan normal form. 

\noindent Let $Y(\theta)=C^*{Z}(\theta)^*J{Z}(\theta)C$. 
Then the coefficients $y_{j,k}(\theta)$ of $Y(\theta)$ 
satisfy $y_{j,k}(\theta)=\langle{z}_j(\theta),Jz_k(\theta)\rangle_\mathbb{C}$ where $z_j(\theta)$ is the $j$-th column 
of $Z(\theta)C$. 
Since $X^t(\theta)^*JX^t(\theta)=J$, then for all $\theta,t$,

\begin{equation}y_{j,k}(\theta)=\langle X^t(\theta){z}_j(\theta),JX^t(\theta)z_k(\theta)\rangle_\mathbb{C}\end{equation}

\noindent Three cases are to be considered:
\begin{enumerate}
\item $y_{j,k}$ is continuous on $\mathbb{T}^d$;
\item $z_j$ is continuous on $2\mathbb{T}^d$ and $z_k$ is continuous on $\mathbb{T}^d$;
\item $z_j$ and $z_k$ are only continuous on $2\mathbb{T}^d$.
\end{enumerate}

\bigskip
Case 1: $z_j$ and $z_k$ are in $\mathcal{W}'$. Then for some $r_j,r_k$,

\begin{equation}\label{gg}\begin{split}&y_{j,k}(\theta)=\langle e^{t(\alpha_j+i\beta_j)}
\sum_{i=r_j}^j\frac{t^{j-i}}{(j-i)!}{z}_i(\theta+t\omega),
Je^{t(\alpha_k+i\beta_k)}\sum_{i=r_k}^k\frac{t^{k-i}}{(k-i)!}z_i(\theta+t\omega)\rangle_\mathbb{C}\\
& =e^{t(\alpha_j+\alpha_k+i\beta_j-i\beta_k)}\sum_{i=r_j}^j\sum_{i'=r_k}^k\frac{t^{j-i}}{(j-i)!}
\frac{t^{k-i'}}{(k-i')!}\langle{z}_i(\theta+t\omega),
Jz_{i'}(\theta+t\omega)\rangle_\mathbb{C}\\
&=e^{t(\alpha_j+i\beta_j+\alpha_k-i\beta_k)}\sum_{i=r_j}^j\sum_{i'=r_k}^k\frac{t^{j+k-i-i'}}{(j-i)!(k-i')!}
y_{i,i'}(\theta+t\omega)\end{split}\end{equation}

\noindent In particular, if $j=r_j$ and $k=r_k$, 

\begin{equation}\label{1d}y_{j,k}(\theta)=e^{t(\alpha_j+i\beta_j+\alpha_k-i\beta_k)}y_{j,k}(\theta+t\omega)\end{equation}

\noindent Developing into Fourier series, since $y_{j,k}$ is continuous on $\mathbb{T}^d$, for all $m\in\mathbb{Z}^d$,

\begin{equation}
\hat{y}_{j,k}(m)=e^{t(\alpha_j+i\beta_j+\alpha_k-i\beta_k)}\hat{y}_{j,k}(m)
e^{2i\pi\langle m,t{\omega}\rangle}\end{equation}

\noindent Thus, either $\hat{y}_{j,k}(m)=0$, or $e^{t(\alpha_j+i\beta_j+\alpha_k-i\beta_k+2i\pi\langle m,\omega\rangle)}=1$ 
for all $t$, 
and then $\alpha_j+i\beta_j+\alpha_k-i\beta_k+2i\pi\langle m,\omega\rangle=0$. But if $m\neq 0$, this is impossible since 
$\beta_j-\beta_k$ is not in $2\pi\langle\mathbb{Z}^d,\omega\rangle\setminus\{0\}$. 
Therefore, $y_{j,k}$ is constant. 

\noindent For any $j,k$, it is possible to show, using equations (\ref{gg}) in the appropriate order, that 
$y_{j,k}$ is constant: equation (\ref{gg}), once developed in Fourier series, gives for all $m\in\mathbb{Z}^d$, 

\begin{equation}\hat{y}_{j,k}(m)=e^{t(\alpha_j+i\beta_j+\alpha_k-i\beta_k)}
\sum_{i=r_j}^j\sum_{i'=r_k}^k\frac{t^{j+k-i-i'}}{(j-i)!(k-i')!}\hat{y}_{i,i'}(m)e^{2i\pi\langle m,t\omega\rangle}
\end{equation}

\noindent Assume $y_{i,i'}$ is constant for all $(i,i')$ such that $i<j$ or $i=j, i'<k$. Then, if $m\neq 0$,

\begin{equation}\hat{y}_{j,k}(m)=e^{t(\alpha_j+i\beta_j+\alpha_k-i\beta_k)}\hat{y}_{j,k}(m)e^{2i\pi\langle m,
t\omega\rangle}
\end{equation}

\noindent which again implies that $y_{j,k}$ is constant.

\bigskip
Case 2: $z_j$ is in $\mathcal{W}$ and $z_k$ in $\mathcal{W}'$. 
Then for some $r_j,r_k$,

\begin{equation}\begin{split}\label{ggg}&y_{j,k}(\theta)
=e^{t(\alpha_j+\alpha_k-i\beta_k)}\sum_{i=r_j}^j\sum_{i'=r_k}^k\frac{t^{j+k-i-i'}}{(j-i)!(k-i')!}
y_{i,i'}(\theta+t\omega)\end{split}\end{equation}

\noindent In particular, if $z_j(\theta)$ and $z_k(\theta)$ generate Jordan subbundles of rank 1, for all 
$\theta,t$,

\begin{equation}\label{wjk}y_{j,k}(\theta)=e^{t(\alpha_j+\alpha_k-i\beta_k)}y_{j,k}(\theta+t\omega)\end{equation}

\noindent Developing this into Fourier series, since $y_{j,k}$ is continuous on $2\mathbb{T}^d$, for all $m\in\mathbb{Z}^d$,

\begin{equation}\label{fou}\hat{y}_{j,k}(m)=e^{t(\alpha_j+\alpha_k-i\beta_k)}\hat{y}_{j,k}(m)
e^{2i\pi\langle m,t\frac{\omega}{2}\rangle}\end{equation}

\noindent So, either $\hat{y}_{j,k}(m)=0$, or $e^{t(\alpha_j+\alpha_k-i\beta_k+i\pi\langle m,{\omega}\rangle)}=1$ for all $t$, 
which implies that $\alpha_j+\alpha_k-i\beta_k+i\pi\langle m,{\omega}\rangle=0$. 
Since $\beta_k$ is not in $\pi\langle\mathbb{Z}^d,\omega\rangle\setminus \{0\}$, this is impossible if $m\neq 0$, so $y_{j,k}$ is constant. 

\noindent For other $j,k$, (\ref{ggg}) implies that $y_{j,k}$ is constant.

\bigskip
Case 3: $z_j$ and $z_k$ are in $\mathcal{W}$. Thus they are in a Jordan basis with real exponent, continuous on 
$2\mathbb{T}^d$.

\noindent If $z_{r_j},\dots z_j$ generate a Jordan subbundle 
with exponent $\alpha_j$ and $z_{r_k},\dots z_k$ generate a Jordan subbundle with exponent 
$\alpha_k$, then for all $\theta,t$, (\ref{ggg}) holds, but with $\beta_k=0$. 

\noindent In particular, if $z_j$ and $z_k$ generate Jordan subbundles of rank 1, for all 
$\theta,t$,

\begin{equation}\label{1d''}y_{j,k}(\theta)=e^{t(\alpha_j+\alpha_k)}y_{j,k}(\theta+t\omega)\end{equation}

\noindent Developing into Fourier series again, since $y_{j,k}$ is continuous on $2\mathbb{T}^d$,

\begin{equation}\hat{y}_{j,k}(m)=e^{t(\alpha_j+\alpha_k)}\hat{y}_{j,k}(m)
e^{2i\pi\langle m,t\frac{\omega}{2}\rangle}\end{equation}

\noindent Thus $y_{j,k}$ is constant.

\noindent More generally, for arbitrary $j,k$, for all $m\in\mathbb{Z}^d$ and all $t$,

\begin{equation}\label{fougen}\begin{split}&\hat{y}_{j,k}(m)=e^{t(\alpha_j+\alpha_k)}
\sum_{i=r_j}^j\sum_{i'=r_k}^k
\frac{t^{j-i}}{(j-i)!}\frac{t^{k-i'}}{(k-i')!}
\hat{y}_{i,i'}(m)e^{2i\pi\langle m,t\omega\rangle}\end{split}\end{equation}

\noindent and we can use these equations in the appropriate order to show that all the coefficients of $Y$ are constant, so $Y$ is 
constant. This implies that ${Z}(\theta)^*J{Z}(\theta)$ does not depend on $\theta$.

\bigskip
$\bullet$
Let $\bar{Z}(\theta)={Z}(\theta){Z}(0)^{-1}$. Then

\begin{equation}\begin{split}&\bar{Z}(\theta)^*J\bar{Z}(\theta)
=({Z}(0)^{-1})^*{Z}(\theta)^*J{Z}(\theta)({Z}(0)^{-1})
=J\end{split}\end{equation} 

\noindent since ${Z}^*J{Z}$ is constant. 
Moreover, $\bar{Z}$ is real, so it is $Sp(2n,\mathbb{R})$-valued. 
It is continuous on $2\mathbb{T}^d$. Finally, for all $\theta,t$,

$$X^t(\theta)={Z}(\theta+t\omega)e^{tA}{Z}(\theta)^{-1}$$

\noindent where $A=\left(\begin{array}{cc}
A_1 & 0\\
0 & A_2\\
\end{array}\right)$ thus 

\begin{equation}X^t(\theta)
=\bar{Z}(\theta+t\omega)e^{t{Z}(0)A{Z}(0)^{-1}}\bar{Z}(\theta)^{-1}\end{equation}

\noindent and therefore, $X$ is $Sp(2n,\mathbb{R})$-reducible modulo 2. $\Box$

\bigskip
This proves theorem \ref{Sp'etc} when $G=Sp(2n,\mathbb{R})$.

\subsection{Orthogonal group}

\begin{prop}\label{ort} Let $X$ be a $GL(n,\mathbb{C})$-reducible cocycle.
If it is $O(n)$-valued, then it is $O(n)$-reducible modulo 2.\end{prop}

\dem It is possible to carry out exactly the same proof as for proposition \ref{CSp}, but defining 
$Y(\theta)$ as $C^*Z(\theta)^*Z(\theta)C$ and not as $C^*Z(\theta)^*JZ(\theta)C$ anymore. This way, its coefficients 
are $y_{j,k}(\theta)=\langle z_j(\theta),z_k(\theta)\rangle_\mathbb{C}
=\langle X^t(\theta)z_j(\theta),X^t(\theta)z_k(\theta)\rangle_\mathbb{C}$; since $X$ is bounded, all the Jordan subbundles 
have rank 
1, 
thus the coefficients $y_{j,k}$ satisfy equations (\ref{1d}), (\ref{wjk}) and (\ref{1d''}) with 
$\alpha_j=\alpha_k=0$. We show in exactly the same way that they are constant, 
then define a function $\bar{Z}$ which is continuous on $2\mathbb{T}^d$ and $O(n)$-valued and such that 
$X^t(\theta)=\bar{Z}(\theta+t\omega)e^{tA}\bar{Z}
(\theta)^{-1}$ for some constant matrix $A$ and for all $t,\theta$. 
$\Box$

\bigskip
This proves theorem \ref{Sp'etc} when $G=O(n)$.

\subsection{$U(n)$-reducibility}

\begin{prop}\label{U} Assume that the continuous cocycle $X$ is $U(n)$-valued and $GL(n,\mathbb{C})$-reducible.
Then $X$ is $U(n)$-reducible.

\end{prop}

\dem By lemma \ref{rédsfJ}, there is a decomposition of $\mathbb{C}^n$ into Jordan subbundles. Since the cocycle 
$X$ is $U(n)$-valued, it is bounded, so all Jordan subbundles have rank 
1 and a purely imaginary exponent. Let $z_1,\dots z_n$ be continuous on $\mathbb{T}^d$, each one generating a Jordan 
subbundle, chosen in such a way that the difference of two exponents cannot be in $2i\pi\langle\mathbb{Z}^d,
\omega\rangle
\setminus\{0\}$. Let $Z(\theta)$ be the matrix whose columns are $z_1(\theta),\dots z_n(\theta)$; 
then there is a diagonal matrix $D$ with coefficients $i\beta_1,\dots i\beta_n$ such that for all $\theta,t$,

$$X^t(\theta)=Z(\theta+t\omega)e^{tD}Z(\theta)^{-1}$$

\noindent Let $Y(\theta)={Z}(\theta)^*{Z}(\theta)$, then the coefficients $y_{j,k}$ of $Y$ satisfy

\begin{equation}y_{j,k}(\theta)=e^{it(\beta_j-\beta_k)}y_{j,k}(\theta+t\omega)\end{equation}

\noindent Developing into Fourier series, for all $n\in\mathbb{Z}^d$,

\begin{equation}\hat{y}(n)_{j,k}=e^{it({\beta_j}-{\beta_k})}\hat{y}(n)_{j,k}\end{equation}

\noindent By construction, ${\beta_j}-{\beta_k}$ is either $0$ or is not in $2\pi\langle\mathbb{Z}^d,\omega\rangle$, 
so $Y$ is constant equal to ${Z}(0)^*{Z}(0)$. Thus, if 
$\bar{Z}(\theta):={Z}(\theta){Z}(0)^{-1}$, then $\bar{Z}(\theta)\in U(n)$ and 

\begin{equation}X^t(\theta)=\bar{Z}(\theta+t\omega)e^{t{Z}(0)D{Z}(0)^{-1}}\bar{Z}(\theta)^{-1}
\end{equation}

\noindent Therefore $X$ is $U(n)$-reducible. $\Box$

\bigskip
This completes the proof of theorem \ref{Sp'etc}.

\section{Discrete cocycles}\label{contdis}

We now want to adapt these results to discrete cocycles. In all this section, we shall assume that $(\omega,1)$ 
is rationally independant.

\bigskip
\Def Let $\bar{\omega}\in\mathbb{R}^D$; $X$ is a continuous (resp. discrete) cocycle over $\bar{\omega}$ 
if it is defined on $\mathbb{T}^D\times\mathbb{R}$ (resp. $\mathbb{T}^D\times\mathbb{Z}$) and for all 
$\bar{\theta}\in\mathbb{T}^D$, $t,s\in \mathbb{R}$ (resp. $t,s\in \mathbb{Z}$), 
$X^{t+s}(\bar{\theta})=X^t(\bar{\theta}+s\bar{\omega})X^s(\bar{\theta})$.

\bigskip
\rem The cocycles we studied in the previous sections are all over $\omega$. But to talk about a discrete cocycle 
over $\omega$, it is necessary to assume that $(\omega,1)$ is rationally independant. 
Notice that if a continuous cocycle $X$ over $(\omega,1)$ is $G$-reducible, then its restriction to integer time 
and to the $d$-dimensional subtorus $\mathcal{T}:=\{(\theta,0),\theta\in\mathbb{T}^d\}$ is a discrete cocycle 
over $\omega$ which is $G$-reducible. Indeed, let $Z:\mathbb{T}^d\rightarrow G$ and $B\in \mathcal{G}$ such that 

\begin{equation}X^t(\theta)=Z(\theta+t(\omega,1))^{-1}e^{tB}Z(\theta)\end{equation}

\noindent It is enough to restrict this expression to integer time and to the subtorus $\mathcal{T}$
to get $G$-reducibility for the discrete cocycle $(n,\theta)\mapsto X^n(\theta,0)$.

\subsection{$G$-exponential discrete cocycles}

Given a discrete $G$-valued cocycle $X$, 
we want to define a suspension of $X$, i.e a continuous $G$-valued cocycle whose restriction to integer times 
and possibly to a subtorus coincides with the initial cocycle. 
But this cannot be done 
if $X$ takes its values in two different connected components of $G$, nor if 
$\theta\mapsto X^1(\theta)$ is not homotopic to the identity in $G$ (since the suspension would be a homotopy). 
However, if there is a $\mathcal{G}$-valued function $A$ which is continuous on 
$\mathbb{T}^d$ such that for all $\theta$, $X^1(\theta)=e^{A(\theta)}$, 
then we can define a continuous $G$-valued cocycle 
whose restriction to integer time and to a subtorus coincides with $X$: this will be done in the following 
proposition.
Recall the definition:

\bigskip
\Def A discrete cocycle $X$ is called $G$-exponential if there exists a $\mathcal{G}$-valued function $A$, continuous on 
$\mathbb{T}^d$, such that $X^1(\theta)=e^{A(\theta)}$ for all $\theta$. 

\begin{prop} Let $X$ be a discrete $G$-exponential cocycle over $\omega$. Then there exists a continuous cocycle 
$\tilde{X}:\mathbb{R}\times \mathbb{T}^d\times\mathbb{T}\rightarrow 
G,(t,\theta,\theta')\mapsto X^t(\theta,\theta')$ over $(\omega,1)$ whose restriction 
to $t\in\mathbb{Z}$ and $\{\theta'=0\}$ coincides with $X$.\end{prop}

\dem By assumption, there exists a $\mathcal{G}$-valued function $A$, continuous on $\mathbb{T}^d$, such that for all 
$\theta$, $X^1(\theta)=e^{A(\theta)}$. 

\noindent For all $(\theta,\theta_{d+1})\in\mathbb{T}^{d}\times [0,1[$, let $B(\theta,\theta_{d+1})=
\phi(\theta,\theta_{d+1})
A(\theta-\theta_{d+1}\omega)$ where  
$\phi$ is a real function continuous on $\mathbb{T}^{d}\times [0,1[$ with support contained in 
$\mathbb{T}^d\times [\frac{1}
{4},\frac{3}{4}]$ such that 

$$\int_0^1\phi(\theta+s\omega,\theta_{d+1}+s)ds=1$$ 

\noindent and for all $n\in\mathbb{Z}$, 
$B(\theta,\theta_{d+1}+n)=B(\theta,\theta_{d+1})$. 
So defined, $B$ is continuous on $\mathbb{T}^d\times \mathbb{R}$ and periodic in $\theta_{d+1}$. 
Let $\bar{B}$ be the continuous function on $\mathbb{T}^{d+1}$ which we obtain by taking the quotient.  

\noindent Let $(t,\theta,\theta_{d+1})\mapsto \tilde{X}^t(\theta,\theta_{d+1})$ be the continuous cocycle satisfying 

$$\frac{d}{dt}\tilde{X}^t(\theta,\theta_{d+1})=\bar{B}(\theta+t\omega,\theta_{d+1}+t)\tilde{X}^t(\theta,\theta_{d+1})$$

\noindent This cocycle is $G$-valued.
Since $\int_0^t\bar{B}(\theta+s\omega,\theta_{d+1}+s)ds$ commutes with $\bar{B}(\theta+t\omega,\theta_{d+1}+t)$ 
for all $\theta,t$, we can compute 
$\tilde{X}^t(\theta,\theta_{d+1})$:

$$\forall t,\theta,\theta_{d+1},\  \tilde{X}^t(\theta,\theta_{d+1})=\exp (\int_0^t \phi(\theta+s\omega,\theta_{d+1}+s)ds A(\theta-\theta_{d+1}\omega))$$

\noindent Thus, for all $\theta\in\mathbb{T}^d$, 

$$\tilde{X}^1(\theta,0)=\exp(A(\theta))=X^1(\theta)$$

\noindent and for $n\in\mathbb{N},n\geq 1$, 

$$\tilde{X}^n(\theta,0)=\tilde{X}^1(\theta+(n-1)\omega,n-1)\dots \tilde{X}^1(\theta,0)=X^1(\theta+(n-1)\omega)\dots X^1(\theta)=X^n(\theta)$$

\noindent and for $n\in \mathbb{Z},n\leq -1$, 

$$\tilde{X}^n(\theta,0)=\tilde{X}^{-n}(\theta+n\omega,n)^{-1}=\tilde{X}^{-n}(\theta+n\omega,0)^{-1}
=X^{-n}(\theta+n\omega)^{-1}=X^n(\theta)$$

\noindent whence the proposition. $\Box$

\bigskip
\rem It is possible to show that if $\theta\mapsto X^1(\theta)$ is homotopic to the identity, 
which is weaker than supposing that $X$ is $G$-exponential, then there exists 
a continuous cocycle whose restriction 
to integer time coincides with $X$. However, this cocycle is not $G$-valued anymore.

\bigskip
\Def The continuous cocycle $\tilde{X}$ defined this way is called a suspension of $X$.

\bigskip
We shall show that $GL(N,\mathbb{C})$-reducibility of a discrete $G$-exponential cocycle implies 
$GL(N,\mathbb{C})$-reducibility of its suspension.

\begin{prop}\label{redsusp} Let $\tilde{X}$ be the suspension of a discrete cocycle $X$ which is $GL(N,\mathbb{C})$-
reducible. Then $\tilde{X}$ is $GL(N,\mathbb{C})$-reducible.\end{prop}

\dem Let $Z\in C^0(\mathbb{T}^d,GL(N,\mathbb{C}))$ and $A\in GL(N,\mathbb{C})$ such that 

$$X^n(\theta)=
Z(\theta+n\omega)^{-1}A^nZ(\theta)$$ 

\noindent for all $\theta,n$. There exists $B\in gl(N,\mathbb{C})$ such that

$$X^n(\theta)=
Z(\theta+n\omega)^{-1}e^{nB}Z(\theta)$$

\noindent Let us define,  
for all $\theta\in\mathbb{T}^{d}$, $\tilde{Z}(\theta,0):=Z(\theta)$, and for all $t\in\mathbb{R}$,

$$\tilde{Z}((\theta,0)+t(\omega,1))=e^{{t}B}\tilde{Z}(\theta,0)\tilde{X}^t(\theta,0)^{-1}$$

\noindent Thus, for all $(\theta,\theta_{d+1})\in\mathbb{T}^{d+1}$,
 
$$\tilde{Z}(\theta,\theta_{d+1})=\tilde{Z}((\theta-{\theta_{d+1}}\omega,0)+{\theta_{d+1}}(\omega,1))
=e^{{\theta_{d+1}}B}Z(\theta-{\theta_{d+1}}\omega,0)
\tilde{X}^{{\theta_{d+1}}}(\theta-{\theta_{d+1}}\omega,0)^{-1}$$

\noindent The map $(\theta,\theta_{d+1})\mapsto \tilde{Z}(\theta,\theta_{d+1})$ is periodic in $\theta$ and for all 
$\theta,\theta_{d+1}$,

\begin{equation}\begin{split}&\tilde{Z}(\theta,\theta_{d+1}+1)= e^{(\theta_{d+1}+1)B}
Z(\theta-{(\theta_{d+1}+1)}\omega,0)
\tilde{X}^{{(\theta_{d+1}+1)}}(\theta-{(\theta_{d+1}+1)}\omega,0)^{-1}\\
& =e^{\theta_{d+1}B}e^{B}Z(\theta-{(\theta_{d+1}+1)}\omega,0)
\tilde{X}^1(\theta-(\theta_{d+1}+1)\omega,0)^{-1}
\tilde{X}^{{\theta_{d+1}}}(\theta-{\theta_{d+1}}\omega,
1)^{-1}\\
& =e^{\theta_{d+1}B}Z(\theta-{\theta_{d+1}}\omega,1)\tilde{X}^{{\theta_{d+1}}}
(\theta-{\theta_{d+1}}\omega,1)^{-1}\\
&=e^{\theta_{d+1}B}Z(\theta-{\theta_{d+1}}\omega,0)\tilde{X}^{{\theta_{d+1}}}
(\theta-{\theta_{d+1}}\omega,0)^{-1}
=\tilde{Z}(\theta,\theta_{d+1})\end{split}\end{equation}

\noindent so $\tilde{Z}$ is periodic in $\theta_{d+1}$. Moreover, for all $\theta,\theta_{d+1},t$,

\begin{equation}\begin{split}&\tilde{X}^t(\theta,\theta_{d+1})=\tilde{X}^{t+{\theta_{d+1}}}(\theta-{\theta_{d+1}}\omega,0)
\tilde{X}^{{\theta_{d+1}}}(\theta-{\theta_{d+1}}\omega,0)^{-1}\\
&=\tilde{Z}((\theta-{\theta_{d+1}}\omega,0)+(t+{\theta_{d+1}})(\omega,1))^{-1}e^{(t+
\theta_{d+1})B}\tilde{Z}(\theta-{\theta_{d+1}}\omega,0)\\
&\tilde{Z}(\theta-{\theta_{d+1}}\omega,0)^{-1}e^{-\theta_{d+1}B}
\tilde{Z}((\theta-{\theta_{d+1}}\omega,0)+{\theta_{d+1}}(\omega,1))\\
&=\tilde{Z}((\theta,\theta_{d+1})+t(\omega,1))^{-1}e^{tB}\tilde{Z}(\theta,\theta_{d+1})\end{split}\end{equation}

\noindent whence the $GL(N,\mathbb{C})$-reducibility of $\tilde{X}$. $\Box$

\bigskip
Now we can form the analogue, for discrete time, if $X$ is a $G$-exponential cocycle, 
of propositions \ref{CR}, \ref{RSL},  
\ref{CSp}, \ref{ort} and \ref{U}. They come as corollaries of the above.

\begin{prop} Let $X$ a discrete $G$-exponential cocycle where $G$ is within $GL(N,\mathbb{R})$,
$SL(N,\mathbb{R})$, $Sp(N,\mathbb{R}), SO(N),SU(N)$ and $GL(N,\mathbb{C})$-reducible. Then 
$X$ is $G$-reducible modulo $\chi_G$, with 

$$\chi_G=\left\{\begin{array}{cc}
2 & \mathrm{if}\ G=GL(N,\mathbb{R}),SL(N,\mathbb{R}),Sp(N,\mathbb{R})\ \mathrm{or}\ G=SO(N)\\
1 & \mathrm{if}\ G=SU(N)\\
\end{array}\right.$$

\end{prop}

\dem Let $\tilde{X}$ be a suspension of $X$. By proposition \ref{redsusp}, $\tilde{X}$ is $GL(N,\mathbb{C})$-
reducible. Moreover, $\tilde{X}$ is $G$-valued, so by propositions \ref{CR}, 
\ref{RSL}, \ref{CSp}, \ref{ort} and \ref{U}, $\tilde{X}$ is $G$-reducible modulo $\chi_G$ with $\chi_G=2$ 
if $G=GL(N,\mathbb{R}), 
SL(N,\mathbb{R}), Sp(N,\mathbb{R})$ or $SO(N)$ and $\chi_G= 1$ if $G=SU(N)$. Thus, $X$ is $G$-reducible modulo 
$\chi_G$. $\Box$

\subsection{General case}

It is possible to extend theorems \ref{CR'} and \ref{Sp'etc} to all discrete cocycles, without even assuming 
that their values are in a connected Lie group, because the proof of the theorems \ref{CR'} and \ref{Sp'etc} 
does not essentially use the fact that time is continuous.  

\bigskip
The definition of a subbundle is the same as in section \ref{sf}. If $X$ is a discrete cocycle, 
an invariant subbundle is a subbundle such that for all $n\in\mathbb{Z}$ and all $\theta\in
\mathbb{T}^d$, 
$X^n(\theta)V(\theta)=V(\theta+n\omega)$. We define a Jordan subbundle of rank $k$ modulo $N$, 
a Jordan subbundle and its exponents in the same way as in section \ref{sf}, but now $t$ varies in 
$\mathbb{Z}$ and not in $\mathbb{R}$ anymore. We show in the same way, using part 2. of lemma \ref{dep}, 
that the exponent of a Jordan subbundle modulo $N$ is well-defined modulo $2i\pi\langle\mathbb{Z}^{d+1},
\frac{(\omega,1)}
{N}\rangle$.
Lemmas \ref{rédsfJ}, \ref{blablabla} and \ref{wj} still hold in the discrete case, but lemma \ref{wj} will be 
reformulated as follows:

\begin{lem} For all $1\leq j\leq r$, there exists $1\leq j'\leq r$ such that $\bar{W_j}=W_{j'}$. 
Moreover, $W_j=\bar{W}_j$ iff $\beta_j\in \pi\langle\mathbb{Z}^d,(\omega,1)\rangle$.
\end{lem}
 
\noindent Proposition \ref{CR} can be reformulated in an analogous way: 

\begin{prop}\label{CRdisc} If $X$ is a real discrete cocycle which is $GL(n,\mathbb{C})$-reducible, 
then there is a decomposition $\mathbb{R}^n=W\oplus W'$ where 
\begin{itemize}

\item $W$ is reducible subbundle modulo 2 with a basis 
$z_1,\dots z_r$ such that for all $(\theta,t)\in\mathbb{T}^d\times \mathbb{Z}$, $X^t(\theta)[z_1(\theta)\ \dots z_r(\theta)]=
[z_1(\theta+t\omega)\ \dots z_r(\theta+t\omega)]e^{tA_1}$ where $A_1$ is a matrix with real spectrum;

\item $W'$ is a reducible subbundle modulo 1 with a basis 
$z_{r+1},\dots z_n$ such that for all $\theta,t$, $X^t(\theta)[z_{r+1}(\theta)\ \dots z_n(\theta)]=
[z_{r+1}(\theta+t\omega)\ \dots z_n(\theta+t\omega)]e^{tA_2}$ where $\sigma(A_2)\cap \mathbb{R}
+i\pi\langle\mathbb{Z}^{d+1}
,(\omega,1)\rangle\setminus\{0\}=\emptyset$ and if $\alpha_1+i\beta_1,\alpha_2+i\beta_2\in \sigma(A_2)$, then 
$\beta_1-\beta_2$ is not in $2\pi\langle\mathbb{Z}^{d+1},(\omega,1)\rangle\setminus \{0\}$. 

\end{itemize}
 \end{prop}

\noindent The proof is exactly the same as in proposition \ref{CR}.

\begin{prop}\label{RSLdisc} If $X$ is a discrete $SL(n,\mathbb{R})$-valued cocycle which is 
$GL(n,\mathbb{R})$-reducible, 
then it is $SL(n,\mathbb{R})$-reducible modulo 2.\end{prop}

\noindent Again, the proof is the same as in proposition \ref{RSL}, except that $t$ varies in 
$\mathbb{Z}$ and not in $\mathbb{R}$ anymore.

\begin{prop}\label{CSpdisc} If $X$ is a discrete $Sp(n,\mathbb{R})$-valued (resp. $O(n)$-valued) cocycle which is 
$GL(n,\mathbb{C})$-reducible (resp. $GL(n,\mathbb{C})$-reducible), then it is 
$Sp(n,\mathbb{R})$-reducible (resp. $O(n)$-reducible) modulo 2. \end{prop}

\noindent The proof is exactly as in propositions \ref{CSp} and \ref{ort}, because the fact that $t$ varies in 
$\mathbb{Z}$ and not in $\mathbb{R}$ does not change the conclusions (we use the second part of lemma 
\ref{dep}).

\begin{prop}\label{udisc} If $X$ is a discrete $U(n)$-valued cocycle which is $GL(n,\mathbb{C})$-reducible, 
then 
$X$ is $U(n)$-reducible.\end{prop}

\noindent The proof is essentially the same as in the continuous case.

\bigskip
Propositions \ref{CRdisc}, \ref{RSLdisc}, \ref{CSpdisc} and \ref{udisc} together give theorem 
\ref{reddisc}.

\section{Applications}

The preceding sections enable us to complete some other results on the full-measure reducibility of a generic 
one-parameter family 
of cocycles.  

\bigskip
\Def $\omega\in\mathbb{R}^d$ is diophantine with constant $\kappa$ and exponent $\tau$, denoted by $\omega\in 
DC(\kappa,\tau)$, 
if for all $n\in\mathbb{Z}^d$, 
$|\langle n,\omega\rangle|>\frac{\kappa}{|n|^\tau}$.

\bigskip
\Def Let $\Lambda$ be an interval of $\mathbb{R}$ and $A\in C^\infty(\Lambda,gl(n,\mathbb{C}))$ a one-parameter family of matrices; 
we say $A$ satisfies the non-degeneracy condition $ND(r,\chi)$ 
on an interval 
$\Lambda$ if there exists $r\in\mathbb{Z}^+$ and $\chi>0$ such that for all 
$\lambda\in\Lambda$, for all $u\in\mathbb{R}$, $\sup_{l\leq r}|\frac{\partial^lg(\lambda,u)}{\partial\lambda^l}|
>\chi$ where 

$$g(\lambda,u)=\prod_{\alpha_i(\lambda),\alpha_j(\lambda)\in \sigma(A(\lambda)),i\neq j}(\alpha_i(\lambda)-\alpha_j(\lambda)
-iu)$$

\Def Let $\Lambda$ an interval of $\mathbb{R}$, denote for $h,\delta >0$ by $\mathcal{F}$ the set of the functions 
defined on $\{z\in\mathbb{C}, |Im z|<h\}\times \{x\in\mathbb{R}, d(x,\Lambda)<\delta\}$, holomorphic 
in the first variable and periodic on the real axis. 


\noindent Let

$$C^\omega_{h,\delta}(\mathbb{T}^d\times \Lambda)=\{f\in\mathcal{F}\ |\ |f|_{h,\delta}:=\sup_{|Im x|< h,d(z,\Lambda)
<\delta}|f(x,z)|<+\infty\}$$

\noindent Finally, let $C^\omega_{h,\delta}(\mathbb{T}^d\times\Lambda, \mathcal{G})$ the set of $\mathcal{G}$-valued maps 
each component of whom is in $C^\omega_{h,\delta}(\mathbb{T}^d\times\Lambda)$.

\subsection{Full-measure reducibility in the symplectic case}

In \cite{HY}, H.He and J.You claim the following: 

\begin{thm}\label{redpp}
Suppose $\omega\in DC(\kappa,\tau)$. Let $A\in C^\infty(\Lambda,gl(n,\mathbb{C}))$ a one-parameter family of matrices satisfying 
the non-degeneracy condition $ND(r,\chi)$ on an interval $\Lambda$. There exists $\epsilon_0>0$ depending on $\kappa$ and $\tau$, 
and there exists $h,\delta$, 
such that if $F\in C^\omega_{h,\delta}(\mathbb{T}^d
\times \Lambda,gl(n,\mathbb{C})),|F|_{h,\delta}\leq \epsilon_0$, then for almost every
 $\lambda\in \Lambda$, the cocycle satisfying

$$\partial_\omega X(\theta)=(A(\lambda)+F(\theta,\lambda))X(\theta)$$

\noindent is $GL(n,\mathbb{C})$-reducible. \end{thm}

\bigskip
Let us also assume that $A(\lambda)\in sp(2n,\mathbb{R})$ for all $\lambda\in\Lambda$ and 
$F\in C^\omega_{h,\delta}(\mathbb{T}^d\times \Lambda,sp(2n,\mathbb{R}))$. Then, as a corollary of proposition 
\ref{CSp} and of H.He and J.You's result, we can reformulate the above in the symplectic case:

\begin{cor}
Suppose $\omega\in DC(\kappa,\tau)$. Let $A(\lambda)$ be a one-parameter family of matrices in $sp(2n,\mathbb{R})$ 
satisfying the non-degeneracy condition $ND(r,\chi)$ on an interval $\Lambda$. There exists $\epsilon_0>0$ depending on 
$\kappa,\tau$, and there exists $h,\delta$, such that if $F\in C^\omega_{h,\delta}(\mathbb{T}^d
\times \Lambda,sp(2n,\mathbb{R})),|F|_{h,\delta}\leq \epsilon_0$, then for almost all $\lambda\in \Lambda$, 
the cocycle satisfying

$$\partial_\omega X(\theta)=(A(\lambda)+F(\theta,\lambda))X(\theta)$$

\noindent is $Sp(2n,\mathbb{R})$-reducible modulo 2. \end{cor}

\subsection{Full-measure reducibility in a compact semi-simple group} 

In \cite{Kr2}, R.Krikorian proved the following theorem:

\bigskip
Suppose $\omega\in DC(\kappa,\tau)$. Let $A$ be a generic element of a compact semi-simple group $G$, $r>0$ and $\Lambda$ 
an interval of $\mathbb{R}$. There exists $\epsilon_0>0$ depending on $\kappa,\tau,\Lambda, A,\omega,r$ such that 
if $F\in C^\omega_r(\mathbb{T}^d,
\mathcal{G})$ and $|F|_r\leq \epsilon_0$, then for almost all $\lambda\in\Lambda$, the cocycle satisfying 

$$\partial_\omega X(\theta)=(\lambda A+F(\theta))X(\theta)$$

\noindent is $G$-reducible modulo an integer $\chi_G$ depending only on $G$. If $G=U(n)$, then $\chi_G=1$.

\bigskip
As a corollary of H.He and J.You's result and of proposition \ref{ort}, we know as well that if $G=O(n)$, 
then $\chi_G=2$. 

\subsection{Does one have full-measure reducibility modulo 1 in any Lie group?}

We first point out the following:

\begin{prop}\label{sansres} If $X$ is a continuous $G$-valued cocycle which is $GL(n,\mathbb{C})$-reducible 
to a cocycle $t\mapsto e^{tB}$ such that the eigenvalues of 
$B$ are not in $\mathbb{R}+i\pi\langle\mathbb{Z}^d,\omega\rangle$, then $X$ is $G$-reducible.
\end{prop}

\dem In the notations of section \ref{redCR}, there is a decomposition of $\mathbb{R}^n$ into 
invariant subbundles $W_1\oplus \dots \oplus W_r$, each $W_j,j\leq r$ being the sum of all Jordan subbundles 
with the same exponent. By assumption on the eigenvalues of $B$, none of the subbundles $W_j$ is its own complex 
conjugate. For all $j$, let $(u_1^j+iv_1^j,\dots u^j_{k_j}+iv^j_{k_j})$ be 
a global basis of $W_j$. Then $(u_1^j,v_1^j,\dots u^j_{k_j},v^j_{k_j})$ is a global basis of 
$(W_j+\bar{W_j})\cap \mathbb{R}^n$. For all $\theta$, let ${Z}(\theta)$ be the matrix whose columns are 
$(u_1^j(\theta),v_1^j(\theta),\dots u_{k_j}^j(\theta),v_{k_j}^j(\theta), 1\leq j\leq r)$,  
then ${Z}$ is continuous on $\mathbb{T}^d$ and $GL(n,\mathbb{R})$-valued and for all $\theta,t$, there exists 
$\tilde{B}$ 
such that 
$X^t(\theta)={Z}(\theta+t\omega)e^{t\tilde{B}}{Z}(\theta)^{-1}$, so $X$ is 
$GL(n,\mathbb{R})$-reducible modulo 1.

\bigskip
If $G=GL(n,\mathbb{R})$, the proof is finished. If $G=SL(n,\mathbb{R}),Sp(n,\mathbb{R})$ or $O(n)$, we do exactly as 
in the proof of \ref{RSL}, \ref{CSp} and \ref{ort}, but since, by assumption, only the case 1 can happen, 
one gets $G$-reducibility modulo 1. $\Box$

\bigskip
\textbf{Question:} Let $A(\lambda)$ be a $\mathcal{G}$-valued one-parameter family 
satisfying a non-degeneracy condition for all $\lambda\in \Lambda$ and $F\in C^\omega_{h,\delta}
(\mathbb{T}^d\times\Lambda)$ sufficiently small. Theorem \ref{redpp} tells that the cocycle $X_\lambda$ 
satisfying 

$$X'_\lambda (t,\theta)=(A(\lambda)+F(\theta,\lambda))X_\lambda (t,\theta)$$

\noindent is $GL(n,\mathbb{C})$-reducible for almost all $\lambda$ to $t\mapsto e^{tB_\lambda}$. 
Is it true that for almost all $\lambda$, the eigenvalues of $B_\lambda$ are not in 
$\mathbb{R}+i\pi\langle\mathbb{Z}^d,\omega\rangle$? If it were the case,
$X_\lambda$ would be $G$-reducible modulo 1 for almost every $\lambda$.

\end{document}